
\input amstex
\documentstyle{amsppt}

\define\ZZ{{\Bbb Z}}
\define\NN{{\Bbb N}}
\define\RR{{\Bbb R}}
\define\QQ{{\Bbb Q}}
\define\CC{{\Bbb C}}

\define\im{\text{Im}~}

\define\rk{\text{rk~}}
\define\M{{\Cal M}}

\define\E{{\Cal E}}
\define\Ha{{\Cal H}}
\define\La{{\Cal L}}


\define\geg{{\goth g}}

\define\tet{{\goth t}}

\magnification1200

\TagsOnRight

\topmatter

\title
Reflection groups in hyperbolic spaces \\
and the denominator formula \\
for Lorentzian Kac--Moody Lie algebras
\endtitle

\author
Viacheslav V. Nikulin \footnote{The author is grateful for
the financial support to
Russian fund of fundamental research; Soros Foundation, grant
MI6000; Soros Foundation and Russian government, grant MI6300;
American Mathematical Society.\hfill\hfill}
\endauthor

\address
Steklov Mathematical Institute,
ul. Vavilova 42, Moscow 117966, GSP-1, Russia.
\endaddress

\email
slava\@nikulin.mian.su
\endemail

\abstract
   This is a continuation of our "Lecture on Kac--Moody Lie algebras of
the arithmetic type" \cite{25}.

   We consider hyperbolic (i.e. signature $(n,1)$) integral symmetric
bilinear form $S:M\times M \to {\Bbb Z}$ (i.e. hyperbolic lattice),
reflection group $W\subset W(S)$, fundamental polyhedron $\Cal M$
of $W$ and an acceptable (corresponding to twisting coefficients)
set $P({\Cal M})\subset M$ of vectors orthogonal to faces of $\Cal M$
(simple roots). One can construct the corresponding Lorentzian Kac--Moody
Lie algebra ${\goth g}=
{\goth g}^{\prime\prime}(A(S,W,P({\Cal M})))$ which is graded by $M$.

   We show that $\goth g$ has good behavior of imaginary roots,
its denominator formula is defined in a natural domain and has good
automorphic properties if and only if  $\goth g$ has so called
{\it restricted arithmetic type}. We show that every finitely generated
(i.e. $P({\Cal M})$ is finite) algebra
${\goth g}^{\prime\prime}(A(S,W_1,P({\Cal M}_1)))$ may be
embedded to  ${\goth g}^{\prime\prime}(A(S,W,P({\Cal M})))$ of the
restricted arithmetic type. Thus, Lorentzian Kac--Moody
Lie algebras of the restricted arithmetic type is a natural class to study.

   Lorentzian Kac--Moody Lie algebras of the restricted arithmetic type
have the best automorphic properties for the denominator function
if they have {\it a lattice Weyl vector $\rho$}. Lorentzian Kac--Moody
Lie algebras of the restricted arithmetic type with generalized
lattice Weyl vector $\rho$ are called {\it elliptic}
(if $S(\rho, \rho)<0$) or {\it parabolic} (if $S(\rho, \rho)=0$).
We use and extend our and Vinberg's results on
reflection groups in hyperbolic spaces to
show that the sets of elliptic and parabolic Kac--Moody Lie algebras
with generalized lattice Weyl vector and lattice Weyl vector are
essentially finite.

    We also consider connection of these results with the recent
results \linebreak by R.Borcherds.
\endabstract

\rightheadtext{Denominator  formula for Kac--Moody algebras}

\leftheadtext{Viacheslav V. Nikulin}

\endtopmatter

\document

\head
\S 0. Introduction
\endhead

In \cite{25}, it was shown that a finitely generated
symmetrizable Kac--Moody algebra has "good behavior" of
imaginary roots if and only if it is finite, affine, of rank two or of
hyperbolic arithmetic type. Hyperbolic arithmetic type means
that Weyl group is a  group generated by reflections in a hyperbolic
space of dimension $\ge 2$ with a fundamental
polyhedron of finite volume.
This paper is a natural continuation of these studies.

We consider a hyperbolic (i.e. of signature $(n,1)$) integral symmetric
bilinear form $S:M\times M \to {\Bbb Z}$ (i.e. a hyperbolic lattice),
a reflection group $W\subset W(S)$, a
fundamental polyhedron $\Cal M$ of $W$ and an acceptable
(corresponding to twisting coefficients) set $P({\Cal M})\subset M$
of vectors orthogonal to faces of $\Cal M$ (simple real roots).
Using these data, one can
construct the corresponding Lorentzian Kac--Moody algebra
${\goth g}={\goth g}^{\prime\prime}(A(S,W,P({\Cal M})))$
which is graded by $M$ (see Sects 2.1 and 2.2).

We show that ${\goth g}={\goth g}^{\prime\prime}(A(S,W,P({\Cal M})))$
has "good behavior" of imaginary roots,
its denominator formula is defined in a natural domain and has good
automorphic properties if and only if this algebra
has so called
{\it restricted arithmetic type} (see Sect. 2.2).
This means that the semi-direct
product $W.A(P(\M ))$ has finite index in $O(S)$. Here
$$
A(P(\M))=\{ g \in O_+(S)\ |\ g(P(\M))=P(\M)\}
$$
is the "group of symmetries" of the fundamental polyhedron.

We show that every finitely generated
(i.e. $P({\Cal M})$ is finite) Lorentzian Kac--Moody algebra
${\goth g}^{\prime\prime}(A(S,W_1,P({\Cal M}_1)))$ has an
embedding to  a Lorentzian
Kac--Moody algebra
${\goth g}^{\prime\prime}(A(S,W,P({\Cal M})))$ of
restricted arithmetic type with the same lattice $S$
(see Sect. 2.2). Thus, it is natural to study Lorentzian Kac--Moody
algebras of restricted arithmetic type.

The denominator formula of a
Lorentzian Kac--Moody algebra of restricted arithmetic type
has the best automorphic properties
if this algebra has a {\it lattice Weyl vector}
$\rho \in M\otimes \QQ$
(see Sect. 2.3 and also 2.4). A lattice Weyl vector is
an element $\rho \in M\otimes \QQ$ such that
$$
S(\rho, \alpha)=-S(\alpha,\alpha )/2\ \text{for any\ }\alpha \in P(\M).
$$
A Lorentzian Kac--Moody algebra with a generalized lattice Weyl vector
$\rho$ (see Definition 1.4.9)
is called {\it elliptic} if it has restricted arithmetic type and
$S(\rho, \rho)<0$; and is called
{\it parabolic} if it has restricted arithmetic type and
$S(\rho, \rho)=0$, and there does not exist
a generalized lattice Weyl vector with negative square.
Ellipticity is equivalent to
finiteness of the index $[O(S):W]$. Parabolicity is
equivalent to restricted arithmetic type and
existence of $0\not=c \in M$ such that $S(c,c)=0$ and
$g(c)=c$ for any $g \in A(\M)$ where $A(\M)$ is infinite.
The corresponding lattice $S$ is called {\it elliptic reflective} and
{\it parabolic reflective} respectively.

We use and extend our's and \'E.B. Vinberg's
results on reflection groups
in hyperbolic spaces to show that sets of primitive elliptic
and parabolic reflective lattices $S$ of $\rk S\ge 3$ are finite
(see Sect. 1.1). For elliptic case this was known 15 years ago. Thus, we
extend this finiteness result for the parabolic case.
Surprisingly, exactly the same method which was used for the
elliptic case is successful for the parabolic one. This shows
that these two cases are very similar, and the method which had
been used for the elliptic case is very natural.

We apply the main geometrical result using to obtain finiteness
results above, to show that the set of elliptic
Lorentzian Kac--Moody algebras
${\goth g}^{\prime\prime}(A(S,W,$ $P({\Cal M})))$ with a
lattice Weyl vector $\rho$ is finite for $\rk S\ge 3$.
For the parabolic case, we obtain the same result if we additionally
suppose that the index $[O(S)_\rho:A(P(\M))]$ $<D$ for a
fixed constant $D>O$. Here $O(S)_\rho$ is the stabilizer subgroup of
$\rho$. See Sect. 1.3. Example 1.3.4 demonstrates
that finiteness may not be true for the parabolic case
without this additional condition.

At last, in Sect. 2.4 we consider connection of these our results
with the recent results by R. Borcherds.

This paper was written during my stay at Steklov Mathematical
Institute, Mos\-cow and Mathematical Institute of G\"ottingen
University (January---March 1995).
I am grateful to these Institutes for hospitality.

I am grateful to Professor \'E.B. Vinberg for very useful remarks.
I am grateful to Professor I. R. Shafarevich for his interest to
and support of these my studies.

I plan to publish this paper in Math. Russian Izvest.

\head
\S 1. Some results on reflection groups of
integral hyperbolic lattices
\endhead

\subhead
1.1. Elliptic and parabolic reflective lattices and reflection groups
\endsubhead

One can consider this section as a complement to our papers
\cite{19}, \cite{20} and to \'E.B. Vinberg \cite{31}.
In \cite{19} and \cite{20} we used signature $(1,n)$ for
a hyperbolic form and used one letter $S$ to denote the form
and the space.  Here we use two letters $S$ and
$M$ to denote a form. Also, here we use signature $(n,1)$ for
a hyperbolic form. These notations are standard for Lie
algebras theory.

\vskip10pt

First, we extend results of the papers above to
"reflection groups with a cusp".

Let
$$
S:M\times M \to \ZZ
\tag1.1.1
$$
be a hyperbolic (i.e. of signature $(n,1)$)
integral symmetric bilinear form over $\ZZ$.
Here $M$ is a free
$\ZZ$-module of a finite rank.
To be shorter, we call $S$ as
{\it hyperbolic lattice}.
We consider the corresponding cone
$$
V(S)=\{ x\in M\otimes \RR \ |\ S(x,x) < 0 \},
\tag1.1.2
$$
choose its half-cone $V^+(S)$ and consider the corresponding
hyperbolic space $\La (S)=V^+(S)/\RR_{++}$ where $\RR_{++}$
denote the set of positive numbers.
Thus, a point of $\La(S)$ is a ray $\RR_{++}x$ where $x\in V^+(S)$.
The distance $\rho $ in $\La(S)$ is defined by the formula:
$$
\cosh \rho (\RR_{++}x, \RR_{++}y)= -S(x,y)/\sqrt{S(x,x)S(y,y)},\
x,y \in V^+(S).
\tag1.1.3
$$
With this distance the curvature is equal to $-1$.
Each element $\delta \in M \otimes \RR $ with
$S(\delta ,\delta) >0$
defines the half-space
$$
\Ha_\delta^+ = \{\RR_{++}x \in \La (S)\
\mid \ S(x, \delta )\le 0\}
\tag1.1.4
$$
which is bounded by the hyperplane
$$
\Ha_\delta = \{ \RR_{++}x \in \La(S)\ \mid \ S(x, \delta )= 0 \}.
\tag1.1.5
$$
The element $\delta \in M\otimes \RR$
is defined by the half-space $\Ha_\delta^+$
(respectively, by the hyperplane $\Ha_\delta$)
up to multiplication on elements of $\RR_{++}$
(respectively, on elements of $\RR^\ast$
of non-zero real numbers).
The $\delta$ is called
{\it orthogonal} to the half-space $\Ha_\delta^+$ (respectively,
to the hyperplane $\Ha_\delta$).

Let $O_+(S)$ be the subgroup of $O(S)$ of the index $2$ which fixes
the half-cone $V^+(S)$. It is well-known (this
an easy corollary of the arithmetic groups theory)
that $O_+(S)$ is discrete in $\La(S)$ and
has a fundamental domain of finite volume. If $\phi \in O_+(S)$
defines a reflection in a hyperplane of $\La(S)$,
then $\phi =s_\delta$ for $\delta \in M$ with $S(\delta , \delta )>0$.
Here
$$
s_\delta (x)=x - (2S(x, \delta)/S(\delta , \delta))\delta,\ \ \ \ x\in S,
\tag1.1.6
$$
and $s_\delta \in O(S)$ if and only if
$$
(2S(M,\delta)/S(\delta ,\delta ))\delta \subset M .
\tag1.1.7
$$
In particular, if
$\delta$ is primitive in $M$, this is equivalent to
$$
S(\delta , \delta) | 2S(M, \delta).
\tag1.1.8
$$
Obviously, $s_\delta$ is the reflection in the hyperplane
$\Ha_\delta$ which is
orthogonal to $\delta$.
The reflection changes places half-spaces $\Ha_\delta^+$ and
$\Ha_{-\delta}^+$. The automorphism $s_\delta\in O(S)$
is called {\it reflection} of the lattice $S$.
Any subgroup of $O(S)$
(or the corresponding discrete group of motions of $\La(S)$)
generated by a set of reflections is called
{\it a reflection group}.

We denote by $W(S)$ the subgroup of $O_+(S)$
generated by all reflections of $S$.
A lattice $S$ is called {\it reflective} if index $[O(S):W(S)]$  is
finite; equivalently, $W(S)$ has a fundamental polyhedron
of finite volume in $\La (S)$.
\'E.B. Vinberg, in particular, proved in \cite{29}
that any arithmetic reflection group $W$ in a hyperbolic
space with the field of definition $\QQ$
is a subgroup $W\subset W(S)$ of finite index for one of
reflective hyperbolic lattices $S$.
For $m \in \QQ$ we denote by $S(m)$ the lattice
which one gets multiplying on $m$ the form of $S$. Obviously,
$S(m)$ is reflective if $S$ does.

Almost 15 years ago there was proved

\proclaim{Theorem 1.1.1 (\cite{19, Appendix, Theorem 1}
and \cite{20, Theorem 5.2.1})}  For a fixed $\rk S\ge 3$,
the set of reflective hyperbolic lattices $S$ is finite up to similarity
$S\ \mapsto \ S(m)$ and isomorphism.
\endproclaim

The proof was based on a purely geometrical result on
convex polyhedra of finite volume in hyperbolic spaces
which we want to formulate.

A convex polyhedron $\M$ in a hyperbolic space $\La(S)$
is an intersection
$$
\M=\bigcap_{\delta \in P(\M)} {\Ha_{\delta}^+}
\tag1.1.9
$$
of several half-spaces orthogonal to elements $\delta \in M\otimes \RR$
with $S(\delta , \delta)>0$. We suppose that $\M$ is locally finite in
$\La(S)$. Then the minimal set $P(\M)$ above is defined canonically up to
multiplication of its elements by positive reals.
Then it is called
{\it the set of vectors orthogonal to faces of $\M$}
(and directed outward) or shortly: {\it orthogonal to} $\M$.
We always suppose that $P(\M)$ has this property.
The polyhedron $\M$ is called
{\it non-degenerate} if it contains a non-empty open
subset of $\La (S)$. A non-degenerate polyhedron $\M$
is called {\it elliptic} (equivalently, it
{\it has finite volume}) if it is a convex envelope of a
finite set of points in $\La (S)$ or at infinity of $\La (S)$.
Then $P(\M)$ is finite.
The proof of Theorem 1.1.1 was based on the following result:

\proclaim{Theorem 1.1.1$^\prime$\ (\cite{19, appendix, Theorem 1})}
Let $\M$ be an elliptic (equivalently, of
finite volume) non-degenerate convex polyhedron in hyperbolic space
$\La (S)$ of $\dim \La (S) =n\ge 2$. Then there are elements
$\delta_1,...,\delta_{n+1}\in P(\M)$ with the following properties:

(a) $\rk [\delta_1,...,\delta_{n+1}]=n+1$;

(b) the Gram diagram of the elements $\delta_1,...,\delta_{n+1}$
is connected (i.e. one cannot divide the set
$\{\delta_1,...,\delta_{n+1} \}$ on two non-empty subsets being
orthogonal to one another).

(c) $-2 \le -2S(\delta_i , \delta_j)/
\sqrt{S(\delta_i,\delta_i)S(\delta_j,\delta_j)}< 62$
(other speaking, we have inequality $-2 \le -S(e_i,e_j)<62$
if we normalize elements $e_i \in P(\M)$ by
the condition $S(e_i,e_i)=2$), $1\le i,j\le n+1$.
\endproclaim

To prove Theorem 1.1.1, one should apply Theorem 1.1.1$^\prime$
to the fundamental
polyhedron $\M$ of $W$ and elements $P(\M)$ which belong to the lattice
$M$. See \cite{20, Theorem 5.2.1}.

Further, we name reflection subgroups $W\subset O(S)$
of finite index and the corresponding reflective hyperbolic lattices $S$
also as
{\it elliptic reflection groups} and
{\it elliptic reflective hyperbolic  lattices} respectively.

\vskip10pt

We want to extend results on elliptic reflection groups and
elliptic reflective lattices above to the following situation.

Let $S$ be a hyperbolic lattice, and $W\subset W(S)$
a reflection subgroup.
Let $\M$ be a fundamental polyhedron of $W$. Let
$$
A(\M) =\{ \phi \in O_+(S)\ |\ \phi (\M)=\M\}
\tag1.1.10
$$
be the group of symmetries of $\M$. Clearly, then the semi-direct
product $W.A(\M) \subset O_+(S)$ where $W$ is a normal subgroup in
$W.A(\M)$.

\definition{Definition 1.1.2} A reflection group
$W\subset W(S)$
is {\it parabolic}
if the group $A(\M)$ is infinite but it
has finite index in $O(S)/W$ (this means that
$W.A(\M)\subset O_+(S)$ has finite index) and there exists
an element $0\not= c\in M$ with $S(c,c)\le 0$ such that
$\phi (c)=c$ for any $\phi \in A(\M)$.
One can easily see (replacing $c$ by $-c$ if
necessary)
that $S(c,c)=0$ and $\RR_{++}c \in \M$ (equivalently,
$S(c, P(\M)) \le 0$). We call the primitive element $c \in M$ (or
the point $\RR_{++}c$ at infinity of $\La (S)$) which satisfies this
condition the {\it cusp} of $W$. One can easily see that the
cusp $c$ (and the point $\RR_{++}c$) is unique.
A hyperbolic lattice $S$ is called {\it parabolic reflective}
if $O_+(S)$ contains a parabolic reflection subgroup $W$.
\enddefinition

We want to prove that Theorem 1.1.1 is also valid for
parabolic reflective hyperbolic lattices $S$.

\proclaim{Theorem 1.1.3}  For a fixed $\rk S\ge 3$,
the set of parabolic reflective hyperbolic lattices $S$
is finite up to similarity $S\ \mapsto \ S(m)$ and isomorphism.
\endproclaim

\demo{Proof} We prove an analog  of Theorem 1.1.1$^\prime$ for
appropriate "parabolic polyhedra". Let us fix a point
$O=\RR_{++}c$ at infinity of $\La (S)$. Thus, $c \in S\otimes \RR$,
$S(c,c)=0$ and $S(c,V^{+}(S))<0$.

Recall that a horosphere
$\E_O$ with the center $O$ is the set of all lines in $\La (S)$
containing the $O$.
The line $l=O\RR_{++}h\in \E_O,\ \RR_{++}h\in \La(S)$,
is the set
$l = \{ \RR_{++}(tc+h)\ \mid \ t\in \RR \ \text{and} \
S(tc+h,tc+h)< 0\} $.
We fix a constant $R>0$. Then there exists a unique $\RR_{++}h\in l$
such that $S(h, c)=-R$ and $S(h,h)=-1$. Let $l_1,l_2\in \E_O$ and
$h_1, h_2$ the corresponding elements
we have defined above. Let
$$
\rho(l_1,l_2)=\sqrt{S(h_1-h_2,h_1-h_2)}.
\tag1.1.11
$$
The horosphere $\E_O$ with this distance is an affine Euclidean space.
If one changes the constant $R$, the distance $\rho$ is multiplying
by a constant. The set
$$
\E_{O,R}=\{ \RR_{++}h \in \La (S)\ \mid
\ S(h,c)=-R\ \text{and} \ S(h,h)=-1 \}
\cup \{ O \}
\tag1.1.12
$$
is a sphere in $\overline{\La(S)}$ touching $\La(S)_\infty$ at the $O$.
Moreover, the sphere $\E_{O,R}$
is orthogonal to every line $l\in \E_O$ at the
corresponding to the $l$ point $\RR_{++}h, h\in \E_{O,R}$. The distance
of $\La(S)$ induces an Euclidean distance in $\E_{O,R}$ which is similar
to the distance \thetag{1.1.11}. The set $\E_{O,R}$ is identified with
$\E_O$ and is also called horosphere.

Let $K\subset \E_O$. The set
$$
C_K=\bigcup_{l\in K}{l}
\tag1.1.13
$$
is called the {\it cone with the vertex $O$ and the base $K$}.

A non-degenerate locally finite polyhedron $\M$ in $\La (S)$ is called
{\it parabolic} (relative to the point $O\in \La(S)_\infty$), if
the conditions 1) and 2) below are valid:

1) $\M$ is finite at the point $O$, i.e. the set
$\{\delta \in P(\M) \mid S(c, \delta) =0\}$ is finite.

2) For any elliptic polyhedron ${\Cal N}\subset \E_O$ (i.e. $\Cal N$ is
a convex envelope of a finite set of points in $\E_O$), the
polyhedron $\M\cap C_{\Cal N}$ is elliptic.

A parabolic polyhedron $\M$ is called {\it restricted  parabolic} if
the set
$$
r (\M) = \{ -S(c, \delta /\sqrt{S(\delta, \delta)}) \ \mid \
\delta\in P( \M) \}
\tag1.1.14
$$
is finite.

Geometrically this means that all
hyperplanes $\Ha_\delta,\ \delta \in P(\M)$, are touching of
a finite set of horospheres with the center $O$.

\proclaim{Theorem 1.1.3$^\prime$} Theorem 1.1.1$^\prime$ is also
valid for any restricted  parabolic
polyhedron $\M$ in hyperbolic space
$\La(S)$ of $\dim \La(S)= n \ge 2$.
Thus, there are elements
$\delta_1,...,\delta_{n+1}\in P(\M)$ with the following properties:

(a) $\rk [\delta_1,...,\delta_{n+1}]=n+1$;

(b) the Gram diagram of the elements $\delta_1,...,\delta_{n+1}$
is connected (i.e. one cannot divide the set
$\{\delta_1,...,\delta_{n+1} \}$ on two non-empty subsets
orthogonal to one another).

(c) $-2\le  -2S(\delta_i , \delta_j)/
\sqrt{S(\delta_i,\delta_i)S(\delta_j,\delta_j)}\le 62$
(other speaking, we have inequality $-2 \le -S(e_i,e_j)\le 62$
if we normalize elements $e_i \in P(\M)$ by
the condition $S(e_i,e_i)=2$), $1 \le i,j \le n+1$.
(Unlike Theorem 1.1.1$^\prime$, one has a non-strict inequality
$\le 62$ to the right.)

A fundamental polyhedron $\M(W)$ of a parabolic reflection group
$W\subset W(S)$ of a hyperbolic lattice $S$
is restricted  parabolic with respect to the cusp $\RR_{++}c$ of
the group $W$, and the above statement holds for $\M(W)$.
\endproclaim

\demo{Proof} For the proof, we normalize orthogonal vectors
$\delta \in P(\M)$ of a polyhedron
$\M$ by the condition $S(\delta, \delta)=2$.
We normalize by this condition orthogonal vectors to all
hyperplanes  and half-spaces below.

To prove Theorem 1.1.1$^\prime$ in \cite{19, Appendix},
we had fixed a point $O$ inside the polyhedron $\M$ and had
used the following formula for
$S(\delta_1,\delta_2)$ of two elements
$\delta_1, \delta_2 \in P(\M)$ with
$S(\delta_1,\delta_1)=S(\delta_2,\delta_2)=2$ (see
\cite{19, Appendix, formula (2.1)}):
$$
-S(\delta_1, \delta_2) =
4{sin{{\theta_1+\theta_{12}}\over 2}sin{\theta_2+\theta_{12}\over2}
\over
{sin{{\theta_1}\over 2}sin{{\theta_2}\over 2}}}-2.
\tag1.1.15
$$
Here $\theta_1$ and $\theta_2$ are angular openings of
cones with a vertex $O$ and bases $\Ha_{\delta_1}$ and $\Ha_{\delta_2}$
respectively. For non-intersecting hyperplanes
$\Ha_{\delta_1}$ and $\Ha_{\delta_2}$, the angle $\theta_{12}$
is the angular opening of the cone with the vertex $O$
which intersects the previous
cones and has minimal angular opening (this is the minimal touching cone of
two cones above). Here the cones and angles are
taken in the sense of hyperbolic geometry. Using analytic continuation,
one can obviously generalize this formula for arbitrary elements $\delta_1,
\delta_2\in M\otimes \RR$
with $S(\delta_1, \delta_1)=S(\delta_2, \delta_2)=2$.
One should consider the plane section containing $O$ and orthogonal to the
hyperplanes $\Ha_{\delta_1}$ and $\Ha_{\delta_2}$,
and define an appropriate orientation of all plane angles of the
section.

To prove Theorem 1.1.3$^\prime$, we use similar formula.
Assume that the polyhedron
$\M$ is restricted  parabolic with respect to the point $O=\RR_{++}c$
at infinity. For an element $\delta \in M\otimes \RR$ with
$S(\delta, \delta)=2$, we consider the
"angle" $\theta (\delta)=-1/S(c , \delta)$.
Then one has the following analog
of the formula \thetag{1.1.15} above.
$$
-S(\delta_1, \delta_2) =
4{{(\theta(\delta_1)+\theta (\delta_{12}))
(\theta (\delta_2)+\theta (\delta_{12}))}
\over
{\theta (\delta_1)\theta (\delta_2)}}-2.
\tag1.1.16
$$
Here like above, for $\delta_1, \delta_2 \in P(\M)$ with
$S(\delta_1, \delta_1)=S(\delta_2,\delta_2)=2$ and
non-intersecting hyperplanes $\Ha_{\delta_1}$, $\Ha_{\delta_2}$, the
element $\delta_{12}\in M\otimes \RR$ is orthogonal to the hyperplane
$\Ha_{\delta_{12}}$ with
$S(\delta_{12},\delta_{12})=2$ and $S(c, \delta_{12})<0$,
which intersects $\Ha_{\delta_1}$ and $\Ha_{\delta_2}$
only by infinite points and has minimal $\theta(\delta_{12})$.

Here $\theta (\delta)$ behaves like an angular opening of the cone
with vertex $O$ and the base $\Ha_{\delta}$. If $A, B, C, O$ are
$4$ consecutive vertices at infinity of a convex polygon on a plane
and $e_1, e_2, e_3$ are orthogonal
to lines $AB, BC$ and $AC$ respectively, have
$S(e_1, e_1)=S(e_2, e_2)=S(e_3, e_3)=2$ and
$-S(c, e_1)>0$, $-S(c, e_2)>0$ and $-S(c, e_3)>0$,
then
$$
\theta (e_3)=\theta(e_1)+\theta(e_2).
\tag1.1.17
$$

Formulae \thetag{1.1.15}, \thetag{1.1.16} and \thetag{1.1.17}
are formulae of
elementary analytic 2-dimensional hyperbolic geometry
if one considers the plane which contains $O$ and is
orthogonal to $\Ha_{\delta_1}$ and $\Ha_{\delta_2}$.

Using the formula \thetag{1.1.1}) and
"angles" $\theta(\delta_1)$, $\theta(\delta_2)$,
$\theta(\delta_{12})$ instead of $\theta_1, \theta_2, \theta_{12}$ for
the formula \thetag{1.1.15}, the
proof of Theorem 1.1.3$^\prime$ for a restricted
parabolic polyhedron $\M$
is completely the same as for Theorem 1.1.1$^\prime$
(see the proof of Theorem 1 in \cite{19, Appendix}).

Let us prove the last statement. Let $K$ be a fundamental polyhedron for
the action of $A(\M)$ on the horosphere $\E_O$.
Let $C_K$ be the cone with the vertex
$O$ and the base $K$. Then $C_K \cap \M(W)$ is a fundamental polyhedron
for the
semidirect product $W.A(\M)$ which has finite index in $O_+(S)$.
Then $C_K\cap \M(W)$ is an elliptic polyhedron. It follows that
$\M(W)$ is a parabolic polyhedron with respect to $O$.
The set of hyperplanes $\Ha_\delta$,
$\delta\in P(\M)$, of faces of the polyhedron $\M(W)$, which are
also hyperplanes of faces of the polyhedron $C_K\cap \M(W)$,
is finite. It follows that $r(\M(W))$ is
finite, and $\M(W)$ is restricted  parabolic.

This finishes the proof of Theorem 1.1.3$^\prime$.
\enddemo

Now Theorem 1.1.3 follows from Theorem 1.1.3$^\prime$
(like Theorem 1.1.1 follows
from Theorem 1.1.1$^\prime$). See the proof of Theorem 5.2.1 in
\cite{20}. In \cite{20} this is done over an arbitrary appropriate field.
Over $\QQ$ the proof is very easy.

This finishes the proof of Theorem 1.1.3.
\enddemo

\`E.B. Vinberg \cite{31} has shown that for an elliptic reflective
hyperbolic lattice
$S$ the rank $\rk S \le 30$. F. Esselmann \cite{8}
improved this result and
has shown that $\rk S \le 22$. This estimate for elliptic reflective
hyperbolic lattices is exact. R. Borcherds \cite{1}
has proved that the maximal even sublattice of the odd unimodular
hyperbolic lattice of the rank $22$ is elliptic reflective.
Using Vinberg's method one can prove that
the rank of parabolic reflective  lattices $S$ is also
absolutely bounded.
An easy estimate one can get is $\rk S \le 43$ for parabolic
reflective lattices $S$.
Here we use existence of the Leech lattice and two different
even unimodular positive lattices of the rank $16$.
But one can expect that the exact estimate here should be
$\rk S \le 26$. J. Conway \cite{7}
proved that the even unimodular hyperbolic lattice of the rank
$26$ is parabolic reflective.

We remark that general results which bound dimension of
arbitrary (not necessarily arithmetic) reflection groups of so
called parabolic and hyperbolic type in hyperbolic spaces were
obtained in \cite{23} (see also \cite{24}). They generalize results of
author \cite{20}, \'E. B. Vinberg \cite{31}, M. N. Prokhorov \cite{27}
and A. G. Khovanskii \cite{14} which bound dimension of
elliptic (i.e. with a fundamental polyhedron of finite volume)
reflection groups in hyperbolic spaces.

\subhead
1.2. Twisting coefficients
\endsubhead

For Kac--Moody algebras and generalized Kac--Moody algebras
which we consider later,  the hyperbolic lattice
$S:M\times M \to \ZZ$ up to similarity
$S \mapsto S(m)$, $m \in \QQ$, is the invariant of the algebra
or of the corresponding generalized Cartan matrix.
(We will consider only indecomposable generalized Cartan matrices.)
Thus, it is natural to normalize $S$ to be primitive or even primitive.

Remind that a lattice $S$ is
{\it even} if $S(x, x)$ is even for any $x \in M$. Otherwise, the lattice
$S$ is called {\it odd}.
The lattice $S$ is {\it primitive} (respectively,
{\it even primitive}) if $S(1/m)$ is not a lattice (respectively
even lattice) for any natural $m \in \NN$ and $m\ge 2$.
Difference between these two normalizations is that if $S$ is primitive
and odd, then $S(2)$ will be primitive even. In most results below, it does
not matter which of these two normalizations is chosen.
Thus, below, "primitive" (lattice)
means primitive or even primitive if we don't
say exactly which normalization we choose.

Let $S:M\times M \to \ZZ$ be a primitive hyperbolic lattice. We fix
a reflection group $W\subset W(S)$ and a fundamental polyhedron $\M$
of $W$. Let $P(\M)_{pr}$ be the set of primitive elements of
$M$ orthogonal to faces of $\M$. Let $A(\M)\subset O_+(S)$
be the group of symmetries of $\M$ (see \thetag{1.1.10}).

\definition{Definition 1.2.1}
A function $\lambda : P(\M)_{pr} \to \NN$ is called
{\it twisting coefficients function} if
$$
S(\lambda (\delta)\delta , \lambda (\delta)\delta )|
2S(M, \lambda (\delta) \delta),\ \text{for any}\ \delta \in P(\M)_{pr}
\tag1.2.1
$$
and the subgroup
$$
A(P(\M ))=\{ g \in A(\M)\ | \lambda (g(\delta))=\lambda (\delta)\
\text{for any}\ \delta \in P(\M)_{pr } \}
\tag1.2.1$^\prime$
$$
has finite index in $A(\M)$.
Number $\lambda (\delta)$ is called
a {\it twisting coefficient} of the $\delta \in P(\M)_{pr}$.
\enddefinition

Clearly, we can reformulate this definition as follows:

Let us consider a new set $P(\M )$ of orthogonal vectors to
$\M$:
$$
P(\M )=\{\alpha =\lambda (\delta)\delta \ | \ \delta \in P(\M )_{pr} \}.
$$
This set $P(\M )\subset M$ is called {\it acceptable}  if
$$
S(\alpha , \alpha )| 2S(M, \alpha )\ \text{for any} \ \alpha \in P(\M )
\tag1.2.2
$$
and the subgroup
$$
A(P(\M ))=\{ g \in A(\M)\ | g(P(\M))=P(\M) \}
\tag1.2.2$^\prime$
$$
has finite index in $A(\M)$.
Obviously, \thetag{1.2.1} and \thetag{1.2.1$^\prime$}
are equivalent to \thetag{1.2.2} and \thetag{1.2.2$^\prime$}.
Thus, an
acceptable sets $P(\M)$ is equivalent to a twisting
coefficient function. If $P(\M )$ is acceptable, then
$$
P(\M)/P(\M)_{pr}
\tag1.2.3
$$
is a twisting coefficient function.

We need to estimate $S(\delta , \delta )$, $\lambda (\delta)$ for
$\delta \in P(\M )_{pr}$, and $S(\alpha, \alpha )$ for $\alpha \in P(\M)$
where $P(\M )$ is acceptable.

\proclaim{Proposition 1.2.2} Let $a(S)$ be the exponent of the discriminant
group $A_S=M^\ast /M$.
Let $\delta \in P(\M )_{pr}$, and $a(\delta )$
be the maximal natural number such that $\delta /a(\delta )\in M^\ast$.
Then $a(\delta)|a(S)$, and the $\lambda (\delta )$ satisfies
\thetag{1.2.1} if and only if
\newline
$\lambda (\delta )S(\delta, \delta )|2a(\delta)$.
In particular,
$\lambda (\delta )S(\delta, \delta )|2a(S)$,
and for any acceptable $P(\M)$ and $\alpha \in P(\M)$ one has
$S(\alpha, \alpha )|4a(S)^2$.
\endproclaim

\demo{Proof} This is trivial.
\enddemo

\remark{Remark 1.2.3}
Let $P(\M )$ be an acceptable set above.
We denote
$$
M^\ast_{P(\M)}=\{ x \in M^\ast\ |\
S(\alpha , \alpha) | 2S(x, \alpha)\ \text{for any}\
\alpha \in P(\M)\}.
\tag1.2.4
$$
We consider elements $h_j\in M^\ast_{P(\M)} \cap \RR_{++}\M$
where $J$ is a countable set.
Then the Gram matrix of elements $P(\M)\cup\{ h_j\ |\ j \in J\}$
defines a generalized Kac--Moody algebra
(see \cite{2}) with the set of simple real roots $P(\M)$, set of
imaginary simple roots $h_j\in M^\ast_{P(\M)},\ j \in J$, and
the Weyl group $W$. One can consider $M$
(respectively,  $M^\ast_{P(\M)}$)
as an "extended root sublattice generated by real simple roots"
(respectively, "extended root lattice"
(generated by real and imaginary simple roots))
modulo kernel of the canonical symmetric bilinear form.
Thus, Proposition 1.2.2 describes all possibilities for the part of this
matrix connected with real simple roots.
\endremark

\subhead
1.3. Elliptic and Parabolic reflection groups with the lattice Weyl vector
\endsubhead

We fix a primitive even elliptic or parabolic reflective lattice
$S:M\times M \to \ZZ$ and elliptic or parabolic reflection
group $W \subset W(S)$. Let $\M$ be a fundamental polyhedron of $W$ and
$P(\M)$ an acceptable set of orthogonal vectors to $\M$. To be
shorter, we name the pair $(S, P(\M))$ {\it elliptic or parabolic
pair} respectively.

\definition{Definition 1.3.1} An element $\rho \in M\otimes \QQ$ is called
{\it lattice Weyl vector} of $P(\M )$ if
$$
S(\rho, \alpha)=-S(\alpha, \alpha)/2\ \text{for any\ } \alpha \in P(\M).
\tag1.3.1
$$
\enddefinition

Evidently, the Weyl vector $\rho \in [P(\M )]^\ast \subset M\otimes \QQ$
if $\rho$ does exist.
Evidently, $\RR_{++}\rho \in \M$  and $\rho $ is invariant with respect to
the automorphism group $A(P(\M))$. It follows that
$S(\rho , \rho )<0$ if $W$ is elliptic, and $\rho \in \QQ_{++}c$ and
$S(\rho , \rho )=0$ if $W$ is parabolic where $c$ is the cusp.
Since $P(\M )$ generates $M \otimes \QQ$ (by Theorems 1.1.1$^\prime$
and 1.1.3$^\prime$),
the Weyl vector $\rho \in M\otimes \QQ$ is evidently
unique for the fixed subset $P(\M)\subset M$.

\vskip10pt

We want to show that the set of elliptic or parabolic
pairs  $(S, P(\M))$ with the lattice Weyl
vector is essentially finite up to isomorphism. Here the {\it pair
$(S, P(\M))$ is isomorphic to a pair
$(S^\prime ,
P^\prime (\M^\prime ))$}
if their exists an isometry $\phi: S \to S^\prime$ of lattices such that
$\phi (P(\M ))=P^\prime (\M^\prime )$.

\proclaim{Theorem 1.3.2} For $\rk S \ge 3$, the set of elliptic pairs
$(S, P(\M))$ with a lattice Weyl vector is finite up to isomorphism.
Here $S$ is a primitive even elliptic reflective lattice,
$W\subset W(S)$ an elliptic (i.e. of finite index in $O(S)$ )
reflection subgroup,
$\M$ a fundamental polyhedron of $W$, and $P(\M)\subset M$ an acceptable
set of all orthogonal vectors to $\M$.
\endproclaim

\demo{Proof} We fix one of primitive even elliptic reflective hyperbolic
lattices $S:M\times M \to \ZZ$ of the $\rk S=n+1\ge 3$
(we already know that their set is finite).
Let $(S, P(\M))$ be an acceptable pair with a lattice Weyl vector.

By Theorem 1.1.1$^\prime$ and
Proposition 1.2.2,  there are
elements $\alpha_1,...,\alpha_{n+1}\in P(\M )$ such that
they generate a sublattice
$M^\prime \subset M$ of finite index and
Gram matrix of these elements has bounded integral coefficients.
Thus, there exists only a finite set of possibilities
for these Gram matrices.

Let us fix one of possibilities above for the Gram matrix
$A= (S(\alpha_i, \alpha_j )),\ \ 1\le i,j \le n+1$. We have:
$M^\prime \subset M\subset (M^\prime)^\ast$ where the embedding
$M^\prime \subset (M^\prime )^\ast$ is defined by the Gram matrix above.
Thus, there exists only a finite set of possibilities for the
intermediate lattice $M$.

Let us fix one of possibilities above for $M$.
Since the lattice $S$ is non-degenerate,
there exists the unique element $\rho \in M\otimes \QQ$
such that $\rho$ satisfies Definition 1.3.1 for the subset
$\{\alpha_1,...,\alpha_{n+1}\}\subset P(\M)$.

By Proposition 1.2.2, $0<S(\alpha, \alpha )\le 4a(S)^2$ is bounded for
$\alpha \in P(\M )$.
It follows that the set
$$
\{\alpha \in M \ | \ 0 > S(\rho, \alpha )=
-S(\alpha, \alpha )/2 \ge -2a(S)^2 \}
\tag1.3.2
$$
is finite. It follows that we have only finite set of possibilities for
the subset $P(\M)\subset M$.

It follows Theorem.
\enddemo

For the parabolic case, the finiteness result will be the following.

\proclaim{Theorem 1.3.3}  For any parabolic pair
$(S, P(\M))$ with a lattice Weyl vector $\rho$ the group
$A(P(\M))\subset O(S)_\rho=\{ \phi \in O(S)\ | \ \phi (\rho)=\rho \}$ has
finite index in the Euclidean crystallographic group $O(S)_\rho$.

We fix a constant $D>0$. Then
for $\rk S \ge 3$, the set of parabolic pairs
$(S, P(\M))$ with a lattice Weyl vector $\rho$ is finite up
to isomorphism if
index $[O(S)_\rho:A(P(\M ))]<D$.
Here $S$ is a primitive even parabolic reflective lattice,
$W\subset W(S)$ a parabolic reflection subgroup with the cusp
$\RR_{++}c =\RR_{++}\rho $, $\M$ a fundamental polyhedron of $W$,
and $P(\M)\subset M$ an acceptable set of orthogonal vectors to $\M$.
\endproclaim

\demo{Proof} Let us prove the first statement.

Let $(S, P(\M))$ be the parabolic pair corresponding to a
parabolic reflection group $W$ with a fundamental polyhedron $\M$ and
an acceptable set $P(\M)$ of vectors
orthogonal to $\M$, and with a Weyl vector
$\rho \in \QQ_{++}c$ where $c$ is the cusp of $W$.

Since $\rho \in \QQ_{++}c$, we have $A(P(\M ))\subset O(S)_c=O(S)_\rho$
and $A(P(\M))$ has finite index in $O(S)/W$. Let $K$ be a fundamental
domain of $A(P(\M))$ on the horosphere $\E_O$ where $O=\RR_{++}c$,
and $C_K$ is the cone with the vertex $O$ and the base $K$.
Then $C_K \cap \M$ is a fundamental domain of finite volume
for $W.A(P(\M ))$ in $\La (S)$. Since
$\rho \in \QQ_{++}c$ is Weyl vector,
then $S(c, \alpha )<0$ for any $\alpha \in P(\M )$.
Thus, there does not exist a face
of $\M$ which contains the cusp $O$; on the other hand,
$O$ belongs to $\M$. It follows
that the fundamental domain $K$ on the horosphere $\E_O$
has finite volume. Since $O(S)_c$ is discrete in $\E_O$, it follows that
$A((P(\M ))$ has finite index in $O(S)_c$.
This proves the first statement.

Let us prove second statement. We are arguing like for the proof of
Theorem 1.3.2 (using Theorems 1.1.3 and 1.1.3$^\prime$
instead of 1.1.1 and 1.1.1$^\prime$).
But at the end of the proof we should use
that number of subgroups $A\subset O(S)_\rho$ is finite if index
$[O(S)_\rho :A] <D$, and
replace the set \thetag{1.3.2} by the set
$$
\{\alpha \in M \ | \ 0>S(\rho, \alpha ) =
- S(\alpha, \alpha )/2 \ge -2a(S)^2
\ \text{and} \  \Ha_{\alpha}\cap C_K \not= \emptyset \}.
\tag1.3.3
$$
Here (like above) $K$ is a fundamental domain for $A(P(\M))$ on the
horosphere $\E_O$, $O=\RR_{++}\rho$, and
$\Ha_{\alpha}$ is the hyperplane which is orthogonal
to $\alpha\in M$.

This finishes the proof.
\enddemo

The next example demonstrates that the condition
$[O(S)_\rho:A(P(\M))]<D$ for
index  is essential in Theorem 1.3.3.

\example{Example 1.3.4}
Let us consider an even primitive hyperbolic lattice
$M=\ZZ\delta_1\oplus \ZZ\delta_2 \oplus \ZZ\delta_3$ of the rank $3$
where
$$
S(\delta_i, \delta_j)=\
\left(
\matrix
2   & -2 & -2 \\
-2  &  2 & -2 \\
-2  & -2 &  2
\endmatrix
\right).
$$
Here $\{ \delta_1, \delta_2, \delta_3 \} = P(\Delta )_{pr}$
where $\Delta$ is a fundamental triangle with zero angles
(i.e, it has $3$ vertices at infinity)
of the reflection group $W^{(2)}(S)$ generated
by reflections in all elements $\delta \in M$ such that
$S(\delta , \delta )=2$. This is one of examples of $2$-reflective
lattices of the rank $3$ which were described in \cite{18}.

Let $O=\Ha_{\delta_2}\cap \Ha_{\delta_3}$. Here $O=\RR_{++}c$ where
$c=\delta_2+\delta_3 \in M$. We have $c^2=0$,
$S(c, \delta_2)=S(c,\delta_3)=0$ and $S(c, \delta_1)=-4$.

Denote
$$
\rho=c/4,\  e_0=\delta_1,\ f_{01}=2s_{\delta_1}(\delta_2),\
f_{02}=2s_{\delta_1}(\delta_3).
$$
One can easily check that
$$
S(\alpha, \alpha )\ |\ 2S(M, \alpha), \ \
S(\rho, \alpha)=-S(\alpha, \alpha )/2
\tag1.3.4
$$
for $\alpha =e_0,\ f_{01},\ f_{02}$. Here $S(e_0,e_0)=2$, and
$S(f_{01},f_{01})=S(f_{02},f_{02})=8$.

Let us consider $\phi \in O(S)_\rho$ such that
$$
\phi (\delta_2)=\delta_3, \phi (\delta_3)=s_{\delta_3}(\delta_2),
\phi (\delta_1)=s_{\delta_3}(\delta_1).
$$
The $\phi$ is the parallel translation on the horosphere $\E_0$ which
sends the line $\Ha_{\delta_2}$ to the line $\Ha_{\delta_3}$ and the
triangle $\Delta$ to the triangle $\phi (\Delta)$ which
has the same vertex $O$ and the same side
$\Ha_{\delta_3}$ with the triangle $\Delta$.

For $k \in \NN$, we define an infinite
fundamental polygon $\M_k$ for a reflection
subgroup $W_k \subset W^{(2)}(S)$ with the acceptable set of vectors
$P(\M_k)$ as follows:
$$
P(\M_k)=\{ \phi^t(e_0)\ | \ t \in \ZZ\ \text{and}\ t\not\equiv 0~\mod k\}
\cup \{ \phi^t(f_{01}), \phi^t(f_{02})\ |\
t\equiv 0~\mod k\}.
\tag1.3.5
$$
Clearly, $P(\M_k)$ is an infinite polygon with zero angles.
It follows that $\M_k$ is a fundamental polygon for a
subgroup $W_k\subset W^{(2)}(S)$ generated by
reflections in all elements
of $P(\M_k )$.
By our construction, $\phi^k\in A(P(\M_k))$
generates a subgroup of finite index in $O_\rho (S)$.
Since \thetag{1.3.4},
the set $P(\M_k)$ is acceptable and $\rho$ is the
Weyl vector for $P(\M_k)$. It follows that $W_k$ is a
parabolic reflection group with the cusp $c$.

Obviously, all pairs $(S, P(\M_k ))$ are different because
elements $\alpha \in P(\M_k)$ have square $S(\alpha, \alpha)=2$ or
$S(\alpha, \alpha)=8$, and exactly
$k-1$ consecutive sides of $\M_k$
have orthogonal vectors from $P(\M_k)$ with the square $2$.
\endexample

\subhead
1.4. Reflection groups of arithmetic type
\endsubhead

We consider hyperbolic lattices $S:M\times M\to \ZZ$ and
reflection groups $W\subset W(S)$.
Let $\M$ be a fundamental polyhedron of $W$ and $P(\M)$
an acceptable set of orthogonal vectors to $\M$.
We want to define a class of these groups which is
interesting from the view-point of Kac--Moody algebras.
{}From the point of view of corresponding Kac--Moody algebras,
the next definition means that imaginary roots
"behave very nice" (see Sect. 2.2 below).

\definition{Definition 1.4.1} Consider an integral cone (semi-group)
$$
Q_+=\sum_{\alpha \in P(\M)} {\ZZ_+\alpha}\subset M
\tag1.4.1
$$
and the corresponding integral dual cone
$$
{Q_+}^\ast=\{ \ x\in M \ | \ S(x, P(\M))\le 0\}\subset M .
\tag1.4.2
$$
The group $W$ has {\it arithmetic type} if
$$
V^+(S)\cap (M\otimes \QQ ) \subset \QQ_{+}Q_+ .
\tag1.4.3
$$
Equivalently, this means that
for any $x \in M$ with $S(x,x)<0$ there exist
$n \in \NN$ and
$a(\delta) \in \ZZ_+$, $\delta \in P(\M),$ which are
almost all equal to $0$ (i.e., only finite
set is not equal to zero), such that
$$
nx =\pm \sum_{\alpha \in P(\M)}{a(\alpha)\alpha}.
\tag1.4.4
$$
Since the cone $\overline{V^+(S)}$ is self-dual, i.e.
$\overline{V^+(S)}^{\ \ast} =\overline{V^+(S)}$, the reflection group
$W\subset W(S)$ has arithmetic type if and only if
$$
{Q_+}^\ast=\{ \ x\in M \ | \ S(x, P(\M))\le 0\}\subset
M\cap \overline{V^+(S)}.
\tag1.4.5
$$
Since the fundamental polyhedron of $W$
$$
\M =(\RR_{++}{Q_+}^\ast \cap V^+(S))/\RR_{++}\subset \La (S)
\tag1.4.6
$$
is locally finitely generated and $M\otimes \QQ$ is everywhere dense in
$M\otimes \RR$, we can reformulate this definitions using real cones.

Consider the cone
$$
\RR_+Q_+=\overline{\sum_{\alpha \in P(\M)}{\RR_+\alpha}}\subset
M \otimes \RR,
\tag1.4.7
$$
and the corresponding dual cone
$$
(\RR_+Q_+)^\ast =\{x \in M\otimes \RR\ | \ S(x, P(\M))\le 0 \}.
\tag1.4.8
$$
Then $W\subset W(S)$ has arithmetic type if and only if
$$
\overline{V^+(S)} \subset \RR_+Q_+,
\tag1.4.9
$$
equivalently,
$$
(\RR_+Q_+)^\ast \subset \overline{V^+(S)}.
\tag1.4.10
$$
Thus, \thetag{1.4.3}, \thetag{1.4.4}, \thetag{1.4.5}, \thetag{1.4.9},
\thetag{1.4.10} are
equivalent definitions for $W$ to have arithmetic type.
Obviously, this definition does not depend on a choice of
an acceptable set $P(\M)$ of orthogonal vectors to $\M$.

One should consider "arithmetic type" as
a very weakened condition of finiteness of volume for
fundamental polyhedron of a reflection group.
Another explanation why reflection groups of
arithmetic type are important is
that they are {\it maximal}: There does not exist a reflection group
$W^\prime\subset W(S)$ with a fundamental polyhedron
$\M^\prime$ such that $W\subset W^\prime$ and
$P(\M)\subset P(\M^\prime)$ and $W\not=W^\prime$ (equivalently,
$P(\M)\not= P(\M^\prime)$). This follows from \thetag{1.4.5}.

\enddefinition

In spite of Definition 1.4.1 is very natural, it
seems that it is too general (gives too "wild" groups in general).
Therefore, we define more narrow class of reflection groups $W$
which have "better automorphic properties" for the denominator
formula of the corresponding Kac--Moody algebras. See
\S 2.

\definition{Definition 1.4.2} A reflection group $W$ has
{\it restricted  arithmetic type} if $W$ has one
reflection at least and  for an acceptable
set $P(\M)$ of orthogonal vectors to $\M$ the group of
symmetries $A(P(\M))$ has finite index in
$O(S)/W(S)$ (more exactly, this means that the corresponding
semi-direct product $W.A(P(\M))$ has finite index in $O(S)$).
\enddefinition

We have the following result which shows that
restricted  arithmetic type implies arithmetic type and gives
many examples of groups $W$ of arithmetic type.

\proclaim{Theorem 1.4.3} Let $S$ be a hyperbolic lattice,
$W \subset W(S)$ reflection group which contains
one reflection at least, $\M$ a fundamental
polyhedron of $W$ and $P(\M)$ an acceptable set of orthogonal
vectors to $\M$ (e.g. one can take $P(\M)_{pr}$).
Assume that the group of symmetries $A(P(\M))$ has finite index in
$O(S)/W$ (i.e. $W.A(P(\M))$ has finite index in
$O(S)$). In other words, $W$ has restricted  arithmetic type.
Then  $W$ has
arithmetic type (i.e. equivalent properties
 \thetag{1.4.3}, \thetag{1.4.4},
\thetag{1.4.5}, \thetag{1.4.9}, \thetag{1.4.10} are valid).
\endproclaim

\demo{Proof} Assume that $W$ does not have arithmetic type.
Then \thetag{1.4.10} is not valid. Thus, there exists an element
$q\in \RR_+Q_+^{\ \ast}$ such that $S(q,q)>0$. We have
$\RR_+\M\subset \RR_+Q_+^{\ \ast}$ where $\M$ is the
fundamental  polyhedron of
$W$ in $\La (S)$ which is locally finite in $\La(S)$ and is non-degenerate.
Here, $\RR_+\M =(\RR_+Q_+^{\ \ast})\cap \overline{V^+(S)}$.
Considering intersection of the convex envelope of
$\RR_+\M$ and $\RR_+q$ with $V^+(S)_\infty$, one can see that
there exists a non-empty open in $\overline{\La(S)}$ subset
$U \subset \M$ such that $U\cap \La(S)_\infty$ is a non-empty open
subset of $\La(S)_\infty$. In particular, for any $\alpha \in P(\M)$
the intersection $\Ha_\alpha\cap U=\emptyset$.

Since $W.A(P(\M))$ has finite index in $O(S)$, a fundamental
domain $\Cal D$ for the action of $A(P(\M))$ in $\M$ is a
fundamental domain for $W.A(P(\M))$ in $\La (S)$. It follows that
$\Cal D$ is a convex polyhedron in $\La(S)$ which is a convex envelope of a
finite set of points in $\La(S)$ and at infinity of $\La(S)$. Since $W$
has one reflection at least,
there exists a codimension one face of
$\Cal D$ which is bounded by a hyperplane $\Ha_{\alpha}$,
$\alpha \in P(\M )$. We consider two cases.

Suppose that $\Cal D$ is bounded.
One can easily see that there exists $g \in A(P(\M))$
such that $g(\Cal D)\subset U$.
Then $\Ha_{g(\alpha)}\cap U$ is non-empty where $g(\alpha)\in P(\M )$.
We get the contradiction.

Now assume that $\Cal D$ has an infinite vertex. It is
well-known that this vertex is  $\RR_{++}c$
where $0\not= c \in M\cap \overline{V^+(S)}$ and $S(c,c)=0$. Besides,
it is known that infinite points
$\RR_{++}g(c)$, $g \in W.A(\M)$, are everywhere dense
in $\La(S)_\infty$. In particular, there exists $g\in A(\M)$
such that $O=\RR_{++}g(c)\in U\cap \La(S)_{\infty}$,
because $U\subset \M$. Then $g({\Cal D})$ is a fundamental domain
of $A(P(\M))$ in $\M$ too. By our construction, hyperplanes of
all codimension one faces of $g({\Cal D})$
which contain $O=\RR_{++}g(c)$, are different from $\Ha_{\alpha}$,
$\alpha \in P(\M)$. It follows that the stabilizer subgroup
$A(P(\M))_{g(c)}$ of $g(c)$ has finite index in $O(S)_{g(c)}$
and has a fundamental domain of finite volume on the horosphere
$\E_O$.
Then hyperplanes $\Ha_{f(g(\alpha))}$ where $f\in A(P(\M))_{g(c)}$, are
everywhere dense in $U\cap \La(S)_\infty$. We again get a contradiction.

This finishes the proof.
\enddemo

\remark{Remark} Theorem 1.4.3 may be considered as
a particular case of some general results about limit sets of
discrete groups in hyperbolic spaces and their normal subgroups.
See L. Greenberg \cite{9} and review \cite{34}.
\endremark

One easily can construct many examples of reflection groups
of restricted arithmetic type and hence, by Theorem 1.4.3,
of arithmetic type. The next construction is
inspired by studying of automorphism groups of Enriques surfaces in
\cite{21}, \cite{22}.

\definition{Definition 1.4.4}
Let us consider a finite symmetric bilinear (or an appropriate quadratic)
form $f:A\times A \to \ZZ/d$ where $A$ is
a finite abelian group and
$d \in \NN$.
A subset $\overline{\Delta} \subset A$ is called a {finite root system in
$(f, A)$} if
there exists an integral lattice $T:N\times N \to \ZZ$, its
sublattice $N_1\subset N$ of finite index, its
reflection subgroup $W\subset W(T)$ and an identification
$$
(N/N_1,\  T~\mod N_1)=(A, f)
\tag1.4.11
$$
such that
$$
\overline{\Delta} =\{\delta~\mod N_1\ |\ s_{\delta} \in W\ \text{and\ }
\delta \ \text{is primitive in }N\} .
\tag1.4.12
$$
It would be interesting to give a reasonable description of finite
root systems.
\enddefinition

Now let us consider a hyperbolic lattice $S:M\times M\to \ZZ$, its
sublattice $M_1\subset M$ of finite index and a finite root
system $\overline{\Delta}$ in $(S~\mod M_1, M/M_1)$. Let
$$
W(\overline{\Delta}) \subset W(S)
\tag1.4.13
$$
be a reflection subgroup generated by all reflections $s_{\delta}$,
where $\delta \in M$ is primitive and
$\delta~\mod M_1 \in \overline{\Delta}$.

Since $\overline {\Delta}$ is finite, using Theorem 1.4.3, we have

\proclaim{Corollary 1.4.5} The reflection subgroup
$W(\overline\Delta)$ has
restricted  arithmetic type (and hence, arithmetic type)
if $W(\overline\Delta )$ has one reflection at least.
\endproclaim

By Proposition 1.2.2, for a reflection $s_\delta$
with a primitive $\delta$ the square
$S(\delta, \delta)$ is bounded. Thus, using Corollary 1.4.5,
for an appropriate finite root system $\overline{\Delta}$,
we formally get the most useful for applications following
statement 1.4.6 below. On the other hand, this statement
follows from Theorem 1.4.3 directly.

\proclaim{Corollary 1.4.6} Let us fix a subset
$\{d_1,...,d_k\}\subset \NN$, and
a subgroup $W\subset W(S)$ is generated by all reflections in
primitive elements $\delta \in M$ such that
$S(\delta, \delta)\in \{d_1,.$ $..,d_k\}$.
Then $W$ has restricted  arithmetic type (and hence,
arithmetic type) if it contains one reflection at least.
\endproclaim

For $S(\delta, \delta)=2$, this statement was proved in \cite{15}.

The next statement shows that any finitely generated
reflection group of a hyperbolic lattice is
contained in a reflection group of
restricted  arithmetic type of the same lattice.
Its interpretation for Kac--Moody algebras see in \S 2.

\proclaim{Theorem 1.4.7} Let $S:M\times M\to \ZZ$ be a hyperbolic
lattice and $W\subset W(S)$ a finitely generated reflection
subgroup with a fundamental polyhedron $\M$ and
finite acceptable set $P(\M)$ of vectors orthogonal to $\M$.
Then $W$ is a reflection subgroup
$W\subset W^\prime\subset W(S)$
of a reflection group $W^\prime$
of restricted  arithmetic type
with a fundamental polyhedron $\M^\prime \subset \M$
and an acceptable set $P(\M^\prime)$ containing $P(\M)$.
\endproclaim

\demo{Proof} Using Proposition 1.2.2, it is sufficient to prove Theorem
for primitive orthogonal vectors $P(\M)_{pr}$ and $P(\M^\prime)_{pr}$.

For a reflection subgroup $\widetilde{W}\subset W(S)$ we denote by
$$
\Delta (\widetilde{W}) =
\{ \text{primitive\ } \delta \in M\ |\ s_{\delta}\in \widetilde{W}\}
\tag1.4.14
$$
the set of all {\it primitive roots} of the reflection
group $\widetilde{W}$.

Let $h \in M$, $\RR_{++}h \in V^+(S)$ and
$0\not\in S(h, \Delta(\widetilde{W}))$.
One has the following Vinberg's algorithm \cite{30}
for calculation of the fundamental polyhedron $\widetilde{\M}$ of
$\widetilde{W}$ (equivalently, the set $P(\widetilde{\M})_{pr})$) which
contains $\RR_{++}h$ (equivalently, $S(h, P(\widetilde\M )_{pr})<0$).
One defines the height
$$
\text{height}~(\delta)=-2S(h, \delta)/\sqrt{S(\delta, \delta)},
\ \delta \in \Delta(\widetilde{W}).
\tag1.4.15
$$
Geometrically, the height is equivalent to the distance between
the point $\RR_{++}h$ and the hyperplane $\Ha_{\delta}$
in the hyperbolic space $\La (S)$. According to \'E.B. Vinberg,
$$
P(\widetilde{\M})_{pr}=\bigcup_{n\in \NN}{P(\widetilde{\M})_{pr}^{\ (n)}}.
\tag1.4.16
$$
where
$$
\split
P(\widetilde{\M})_{pr}^{\ (n)}=&\{ \delta\in \Delta(\widetilde{W})\ |\
\text{height}(\delta)=\sqrt{n}\\
 &\text{and}\ S(\delta, P(\M)_{pr}^{\ (m)})\le 0\
\text{for all\ } 1\le m  <n \} .
\endsplit
\tag1.4.17
$$
It is true also that the set
$$
\{\delta \in \Delta (W(S))\ |\ 0\le \text{height}~(\delta)\le C \}
\tag1.4.18
$$
is finite for the fixed $C>0$.

Now let us prove Theorem. We fix $h\in M$ such that
$\RR_{++}h \in \M$ and $S(h, $ $P(\M)_{pr}) <0$ (i.e.
$\RR_{++}h$ is inside $\M$). We use this $h$ to define the height.
Let $N$ be the largest height of elements from the finite
set $P(\M)_{pr}$.

We consider a sublattice $M_1\subset M$ of finite index and the finite
symmetric bilinear form
$$
(S \mod M_1, M/M_1)
\tag1.4.19
$$
with the finite root system
$$
\overline{\Delta(W)}=\Delta(W)~\mod~M_1
\tag1.4.20
$$
such that the set
$$
\{\delta \in \Delta (W(S))\ |\
0\le \text{height}~(\delta)\le N
\ \text{and\ } \delta~\mod M_1 \in \overline{\Delta(W)}
\}
$$
is equal to the set
$$
\{\delta \in \Delta (W)\ |\
0\le \text{height}~(\delta)\le N\}.
$$
Since the set \thetag{1.4.18} is finite,
we always can satisfy this equality for a
sufficiently deep sublattice $M_1\subset M$.
Let us apply Vinberg's algorithm for calculation of $P(\M)_{pr}$ and
the fundamental polyhedron $\M_1$ (equivalently,
$P(\M_1)_{pr}$ ) of the reflection group $W(\overline{\Delta(W)})$.
Since the equality above, we evidently have that
$P(\M_1)_{pr}^{\ (n)}=P(\M)_{pr}^{\ (n)}$
for all $0\le n \le N$. It follows that
$P(\M)_{pr}\subset P(\M_1)_{pr}$. By Corollary 1.4.5, the reflection
group $W(\overline{\Delta(W)})$
has restricted  arithmetic type because
it contains the non-trivial reflection subgroup $W$.

This finishes the proof.
\enddemo

Using the same idea, one can prove the following more general statement:

\proclaim{Theorem 1.4.8} Let $S:M\times M\to \ZZ$ be a hyperbolic lattice
and $W\subset W(S)$ a finitely generated reflection
subgroup with a fundamental polyhedron $\M$ and
finite acceptable set $P(\M)$ of vectors orthogonal to $\M$. Let
$h_1,...,h_k \in (M-\{0\})\cap \RR_{++}\M$ (equivalently,
$h_i \in (M-\{0\})\cap \overline{V^+(S)}$ and
$S(h_i, P(\M))\le 0$, $i=1,...,k$).
Then $W$ is a reflection subgroup
$W\subset W^\prime\subset W(S)$
of a reflection group $W^\prime$
of restricted  arithmetic type
with a fundamental polyhedron $\M^\prime \subset \M$
and an acceptable set $P(\M^\prime)$ which contains $P(\M)$. Moreover,
$h_1,...,h_k \in (M-\{0\})\cap \RR_{++}\M^\prime$.
\endproclaim

\vskip10pt

At the end of this section, we consider relation between elliptic and
parabolic reflection groups and reflection groups of arithmetic
and restricted arithmetic type.

We have the following natural generalization of Definition 1.3.1.

\definition{Definition 1.4.9}
Let $S:M\times M\to \ZZ$ be a hyperbolic lattice
and $W\subset W(S)$  a reflection
group with a fundamental polyhedron $\M$ and
an acceptable set $P(\M)$ of vectors orthogonal to $\M$.

An element $0\not= \rho \in M$ is called a
{\it generalized lattice Weyl vector} of $P(\M )$ if there exists
$C>0$ such that
$$
0 \le -S(\rho, \alpha) \le C \ \text{for any\ } \alpha \in P(\M).
$$

An element $\rho \in M\otimes \QQ$ is called a
{\it lattice Weyl vector} of $P(\M )$ if
$$
S(\rho, \alpha)=-S(\alpha, \alpha)/2\ \text{for any\ } \alpha \in P(\M).
\tag1.4.21
$$

Clearly, some multiple of a lattice Weyl vector gives a generalized
lattice Weyl vector. Clearly, existence of a generalized lattice
Weyl vector does not depend from a choice of an acceptable set
$P(\M)$ of vectors orthogonal to $\M$.
\enddefinition

We have the following statements which mean that elliptic
reflection groups are exactly reflection groups of
arithmetic type with a generalized lattice Weyl vector $\rho$
which has negative square $S(\rho , \rho)$.
And parabolic reflection groups
are exactly reflection groups which have restricted arithmetic type
and have a generalized lattice Weyl vector with zero square
$S(\rho, \rho )=0$
and do not have a generalized lattice Weyl vector with negative square.

\proclaim{Proposition 1.4.10} Any elliptic or parabolic reflection
group $W \subset W(S)$ of a hyperbolic lattice $S$
has restricted arithmetic type
(and by Theorem 1.4.3, it has arithmetic type). If $W$ is
elliptic, it has a generalized lattice Weyl vector $\rho$ with
$S(\rho , \rho )<0$. If $W$ is parabolic, it has a generalized
lattice Weyl vector $\rho $ with $S(\rho , \rho )=0$ and does not have
a generalized lattice Weyl vector with negative square.

\endproclaim

\demo{Proof} This group has restricted
arithmetic type by definition.
By Theorem 1.4.3, it has arithmetic type. Let $\M$ be
a fundamental polyhedron of $W$. If $W$ is elliptic,
every element
$0\not= \rho \in M$ such that
$\RR_{++}\rho \in \text{Int\ }\M$ (equivalently,
$S(\rho, P(\M)) < 0$) is a generalized lattice
Weyl vector. If $W$ is parabolic, the cusp $c$ is a generalized
lattice Weyl vector with $S(c,c)=0$. Below we
prove that $W$ does not have a generalized lattice
Weyl vector with negative square.
\enddemo

\proclaim{Proposition 1.4.11}
Let $S$ be a hyperbolic lattice
and $W\subset W(S)$  a reflection group of arithmetic type,
$\M$ a fundamental polyhedron of $W$, and $P(\M)$ an acceptable
set of orthogonal vectors to $\M$.
Assume that $P(\M)$ has a generalized lattice
Weyl vector $\rho $.

Then $S(\rho , \rho )\le 0$.
If additionally $S(\rho , \rho )<0$, then $W$ is elliptic;
if $S(\rho , \rho )=0$ and $W$ has restricted arithmetic type,
then $W$ is either elliptic or parabolic.
\endproclaim

\demo{Proof} By definition of generalized lattice Weyl vector $\rho $,
we have
$S(\rho, P(\M))\le 0$. Thus, $\rho \in {Q_+}^\ast$ (see Definition 1.4.1).
By \thetag{1.4.5}, $\rho \in \overline{V^+(S)}$. It follows that
$S(\rho , \rho )\le 0$ and $\RR_{++}\rho \in \M$.

Assume that $S(\rho , \rho )<0$. By definition of generalized lattice
Weyl vector, it then follows that $P(\M)$ is finite
(since the lattice $S$ is hyperbolic). By \thetag{1.4.6} and
\thetag{1.4.10}, the fundamental polyhedron $\M=(\RR_+Q_+)^\ast/\RR_{++}$
is a convex envelope of a finite set of points of $\La(S)$ and
$\La(S)_\infty$. Thus, $\M$ has finite volume and $W$ is elliptic.

Now assume that $S(\rho , \rho )=0$ and
$W$ has restricted arithmetic type.
If $A(\M)$ is finite, then $W\subset O(S)$ has finite index
and $W$ is elliptic. Assume that $A(\M)$ is infinite. Obviously,
$g(\rho )$ is also a generalized lattice Weyl vector for any $g\in A(\M)$.
If $g(\rho )\not=\rho $, it follows that $P(\M)$ is finite. Since
$W$ has arithmetic type, like above, one can prove that $W$ is
elliptic and $A(\M)$ is finite. We get a contradiction. Thus, we
have proven that $g(\rho )=\rho $ for any $g\in A(\M)$. Thus, $\rho $ is
the cusp for $A(\M)$. Since $W$ has restricted arithmetic type,
$W.A(\M)$ has finite index in $O(S)$. It follows that $W$ is parabolic.
\enddemo

\remark{Remark 1.4.12} Modifying Example 1.3.4, one can construct
an example of a reflection group $W$ of arithmetic type and with
a lattice Weyl vector $\rho $ with $S(\rho ,\rho )=0$ such that this group
does not have restricted arithmetic type and is not then parabolic.
\endremark

\head
\S 2. Denominator formula for Lorentzian Kac--Moody algebras
\endhead

We restrict considering Kac--Moody algebras, but similarly one can
consider generalized Kac--Moody algebras (see Remark 1.2.3 and
Theorem 1.4.8) and, it seems, appropriate
generalized Kac--Moody superalgebras.

\subhead
2.1. Basic definitions on Lorentzian Kac--Moody algebras
\endsubhead

We refer to V. Kac \cite{13} and R. Borcherds \cite{2} for
basic definitions on Kac--Moody algebras.

We fix an even primitive hyperbolic lattice $S:M\times M \to \ZZ$,
a finitely generated reflection group $W\subset W(S)$,
its fundamental polyhedron $\M$ and an acceptable set $P(\M)$
of vectors orthogonal to $\M$.  Since $W$ is finitely generated,
the set $P(\M)$ is finite. We assume that $P(\M)$ generates a sublattice
of $M$ of finite index and Gram graph of $P(\M)$ is connected.
The last means that one
cannot divide $P(\M)$ on two orthogonal non-empty
subsets. We require this for not to consider finite-dimensional
and affine Lie algebras.

Let us consider a finite matrix
$$
A=(a_{ij})=\left( {2S(\alpha_i,\alpha_j)\over
S(\alpha_i,\alpha_i)} \right),\
\alpha_i, \alpha_j \in P(\M).
\tag2.1.1
$$
Since $\M$ is a fundamental polyhedron of a reflection group and
$P(\M)$ is acceptable, one has that
$a_{ii}=2$, $a_{ij}\le 0$ are integral for $i\not=j$, and
$a_{ij}, a_{ji}$ are equal to $0$ simultaneously.
Matrices $(a_{ij})$ with these properties are called
{\it generalized Cartan matrices}. The generalized Cartan matrix
$(a_{ij})$ defined in \thetag{2.1.1} has special type.
First, it is {\it symmetrizable} which means that
$$
A=DB
\tag2.1.2
$$
where $D$ is diagonal with positive rational coefficients
and $B$ is symmetric. Evidently,
for our case,
$$
D=diag~(..., 2/S(\alpha_i,\alpha_i),...)
\tag2.1.3
$$
and
$$
B=(b_{ij})=\left(S(\alpha_i,\alpha_j)\right)
\tag2.1.4
$$
is the Gram matrix of $P(\M)$. Since the Gram graph defined by $B$
is connected, the matrix $A$ is {\it indecomposable}.
Besides, for our case, the matrix $B$ is {\it hyperbolic} which
means that it has exactly one negative square (and several positive and
zero squares).

A generalized Cartan matrices with properties above is called
{\it indecomposable Lorentzian symmetrizable}. One can easily see
that by the construction above,
any indecomposable Lorentzian symmetrizable generalized
Cartan matrix $A$ defines the primitive even hyperbolic lattice $S$,
$W\subset W(S)$ and $P(\M)$ canonically up to isomorphism and
a finite set of possibilities.
If one requires that $P(\M)$ generates $M$,
then $S,W, P(\M)$ are defined uniquely up to isomorphism.

Let $P(\M)=\{ \alpha_1,...,\alpha_N\}$.
The generalized Cartan matrix $A$ above defines
the Kac--Moody algebra
$\geg^\prime (A)$ over $\CC$.
Lie algebra $\geg^\prime (A)$ is defined by
$3N$ generators $e_1,...,e_N$, $f_1,....,f_N$, $h_1,...,h_N$ and
defining relations
$$
[h_i,h_j]=0,\ [e_i,f_j]=\delta_{ij}h_j,\ [h_i,e_j]=a_{ij}e_j,
[h_i,f_j]=-a_{ij}f_j,
$$
$$
(\text{ad\ } e_i)^{1-a_{ij}}e_j=0,\ (\text{ad\ } f_i)^{1-a_{ij}}f_j=0,\
\text{if $i\not=j$},
\tag2.1.5
$$
for all $1\le i,j \le N$.

The basic property of this Kac--Moody algebra is that
$$
\geg^{\prime\prime} (A) =
\geg^\prime (A)/\tet =\bigoplus_{\alpha \in M}\geg_\alpha
\tag2.1.6
$$
is simple and graded by the lattice $M$. Here $\tet$ is
the center of $\geg^\prime (A)$; $\dim \geg_\alpha <\infty$,
$[\geg_\alpha, \geg_\beta]\subset \geg_{\alpha+\beta}$;
$\geg_0=\CC h_1+...+\CC h_N=M^\ast\otimes \CC$ where
$<h_i,\alpha_j>=a_{ij}$; an element $X\in \geg_\alpha$ iff
$[h,X]=<h,\alpha>X$, for any $h \in \geg_0$. The form $S$
may be extended canonically
to $\geg^{\prime\prime}(A)$ to be an invariant symmetric
bilinear form. It is called {\it canonical symmetric bilinear form}.

The part $\geg_0$ is called {\it Cartan subalgebra} of
$\geg^{\prime\prime}(A)$.
An element $0\not=\alpha \in \M$ is called {\it root} if
the multiplicity $\text{mult}(\alpha)=\dim \geg_\alpha >0$.
Let $\Delta$ be the set of all roots. Roots and their multiplicities
are invariant with respect to $W$ which is called {\it Weyl group}.

One has the following description of the set of roots
$\Delta$ according to V. Kac \cite{13}.
A root $\alpha\in M$ is called {\it real} if $S(\alpha, \alpha )>0$.
A root $\alpha \in M$ is called {\it imaginary}
if $S(\alpha, \alpha )\le 0$.
The set of real roots is denoted by $\Delta^{re}$,
the set of imaginary
roots is denoted by $\Delta^{im}$.
Let $Q_+=\ZZ_+\alpha_1+\cdots +\ZZ_+\alpha_N$ be
an integral cone generated by
simple real roots (we have defined and used this cone in
\thetag{1.4.1}).
We have $\Delta =\Delta_+ \cup \Delta_-$ where
$\Delta_+=\Delta\cap Q_+$, $\Delta_-=-\Delta_+$. Similarly one defines
$\Delta^{re}_+, \Delta^{re}_-$ and $\Delta^{im}_+, \Delta^{im}_-$.
Obviously, a root $\alpha\in \Delta_+$ if and only if
$S(\alpha, h)\le 0$
for any $\RR^+h \in \M$. An imaginary root $\alpha \in \Delta_+$ if
and only if
$\alpha \in \overline{V^+(S)}$.
Elements $\{\alpha_1,...,\alpha_N\}=P(\M )$ are real roots.
They are called {\it simple real
roots} and have multiplicity one.  We have
$$
\Delta^{re}=W(P(\M)).
\tag2.1.7
$$
Let
$$
{Q_+}^\ast=\{ x \in M\ | \ S(x,P(\M))\le 0\}
\tag2.1.8
$$
be the integral cone which is dual to $Q_+$, and
$$
K=Q_+\cap {Q_+}^\ast.
\tag2.1.9
$$
By V. Kac \cite{13} (see also considerations in \cite{25}), one has
$$
\Delta^{im}_+=W(K).
\tag2.1.10
$$
Let us consider the corresponding real cones
$$
\RR_+Q_+,\  (\RR_+Q_+)^\ast\  \text{and}\
\RR_+K=(\RR_+Q_+)\cap (\RR_+Q_+)^\ast.
\tag2.1.11
$$
We have
$$
\RR_+K \subset \RR_+\M .
\tag2.1.12
$$
The cone
$$
T^\ast=W(\RR_+K)=W(\RR_+Q_+\cap (\RR_+Q_+)^\ast)\subset \overline{V^+(S)}
\tag2.1.13
$$
is called {\it dual Tits cone}. The cone $T=(T^\ast)^\ast$ is
called {\it Tits cone}. Obviously,
$$
\Delta^{im}_+=Q_+\cap T^\ast.
\tag2.1.14
$$

We have the following Weyl--Kac denominator formula:
$$
\Phi (z)=:
\sum_{w \in W} {\text{det\ }(w)e^{-2\pi i S(w(\rho ), z)}}=
e^{-2\pi iS(\rho,z)} \prod_{\alpha \in \Delta_+}
{(1-e^{-2\pi i S(\alpha , z)})^{\text{mult}~(\alpha)}}.
\tag2.1.15
$$
Here $\rho \in (\oplus \ZZ \alpha_i)^\ast$  is any element which
satisfies the equality
$$
<\rho, \alpha_i>=-S(\alpha_i, \alpha_i)/2\
\text{for any\ } \alpha_i \in P(\M).
\tag2.1.16
$$
It is called {\it Weyl vector}.
The variable
$$
z \in M\otimes \CC.
\tag2.1.17
$$
The domain  $\Omega (S)=M\otimes \RR + iV^+(S)\subset M\otimes \CC$
is called {\it the complexified cone $V^+(S)$}.
It is known (by V. Kac \cite{13}) that $\Phi (z)$ converges absolutely
in the complexified Tits cone
$$
\Omega (S, W)= M\otimes \RR +i\text{Int\ } T\supset \Omega(S)
\tag2.1.18
$$
if $\text{im\ }z>>0$, and diverges absolutely in
$S\otimes \CC-\overline{\Omega (S,W)}$.
Thus, the domain $\Omega (S, W)$ where the function $\Phi (z)$
converges always contains the natural domain $\Omega (S)$.

We denote the generalized Cartan matrix
$A$ in \thetag{2.1.1} as $A(S,W,P(\M))$ and
corresponding Kac--Moody algebras $\geg^\prime(A)$
and $\geg^{\prime\prime}(A)$ as $\geg^\prime (A(S,W,P(\M)))$
and $\geg^{\prime\prime}(A(S,W,P(\M)))$. They are also called
Lorentzian Kac--Moody algebras.

\vskip10pt

One can see that Lorentzian
Kac--Moody algebras
$\geg^{\prime\prime}(A(S,W,P(\M)))$
have bad properties in general.

\vskip10pt

(-a) Rational cone $\QQ_+\Delta^{im}_+$ generated by imaginary roots gives
only a part of the natural cone $(M\otimes \QQ)\cap V^+(S)$.

\vskip5pt

(-b) The denominator function $\Phi (z)$ converges in the bigger domain
$\Omega (S,W)$ than the natural domain $\Omega(S)$.

\vskip5pt

(-c) The denominator function $\Phi (z)$ is anti-invariant with respect to
the Weyl group $W$ and is invariant with respect to a finite group
$A(\M)$ of symmetries of $\M$. The corresponding semi-direct
product $W.A(S)$ may be very small in the natural arithmetic group
$O(S)$.

\vskip5pt

(-d) The expressions $S(\rho, z)$ and $S(w(\rho), z)$ are not correctly
defined if there does not exist
$\rho \in M\otimes \QQ$ such that
$$
S(\rho, \alpha_i)=-S(\alpha_i, \alpha_i)/2.
$$
If such a $\rho \in M\otimes \QQ$ does exist, we say that $S$ {\it has
a lattice  Weyl vector}. If a lattice Weyl vector does not exist,
one should rewrite \thetag{2.1.15} as
$$
\sum_{w \in W} {\text{det\ }(w)e^{-2\pi i S(w(\rho)-\rho, z)}}=
\prod_{\alpha \in \Delta_+}
{(1-e^{-2\pi i S(\alpha , z)})^{\text{mult}~(\alpha)}}.
\tag2.1.19
$$
Then it is defined correctly.

\vskip10pt

The only way to improve bad properties (-a), (-b) and (-c)
is to consider bigger reflection groups
$$
W=W_1\subset W_2\subset \cdots \subset W_k\subset \cdots \subset W(S)
\tag2.1.20
$$
with fundamental polyhedra $\M_k$ and acceptable sets
$P(\M_k)$ such that
$$
P(\M)=P(\M_1)\subset P(\M_2)\subset \cdots \subset P(\M_k)\subset \cdots .
\tag2.1.21
$$
These gives an increasing sequence of Lorentzian Kac--Moody algebras
$$
\geg^{\prime\prime}(A(S, W_1, P(\M_1)))\subset
\geg^{\prime\prime}(A(S, W_2, P(\M_2)))\subset  \hskip-3pt \cdots \subset
\geg^{\prime\prime}(A(S, W_k, P(\M_k)))\subset \hskip-3pt \cdots
\tag2.1.22
$$
with the same Cartan subalgebra $M^\ast \otimes \CC$. One can
see that this procedure gives an increasing sequence of
dual Tits cones
$$
T_1^\ast \subset T_2^\ast \subset \cdots \subset T_k^\ast \subset
\cdots \subset \overline{V^+(S)}
\tag2.1.23
$$
and then a decreasing sequence of Tits cones $T_k$
and domains $\Omega(S, W_k)$.
Thus, at least properties (-a) and (-b) will be improving.
Unfortunately, this procedure is infinite in general
and gives an infinite generated Lorentzian Kac--Moody algebra.

In the next Section, we will show that using this may be infinite
procedure we can always reverse bad properties (-a), (-b) and (-c).

\subhead
2.2. Lorentzian Kac--Moody algebras of arithmetic and
restricted arithmetic type
\endsubhead

Let $S:M\times M\to \ZZ$ be a hyperbolic lattice, $W\subset W(S)$
an arbitrary reflection subgroup, $\M$ a fundamental polyhedron
of $W$ and $P(\M)$ an acceptable set of
orthogonal vectors to $\M$. Like above, we suppose that $P(\M)$
generates a sublattice of finite index in $M$ and has connected
Gram graph.

Using the same definition as in Sect. 2.1, one can define
an indecomposable hyperbolic symmetrizable generalized
Cartan matrix $A=A(S,W,P(\M))$, and Lorent\-zian Kac--Moody
algebras
$$
\geg^\prime (A)=
\geg^\prime(A(S,W,P(\M))),\ \
\geg^{\prime\prime}(A)=\geg^{\prime\prime}(A(S,W,P(\M))).
$$
They have the same properties as in Sect. 2.1. The only difference
is that the set $P(\M)$ is infinite, the matrix $A$ is infinite and
the set of generators $e_k, f_k, h_k$, where $k \in P(\M)$, is infinite.
This Lie algebra is a union of
its Kac--Moody subalgebras which are defined by finite sets
$P(\M)^{(\le n)}\subset P(\M)$ such that their elements
have height $\le \sqrt{n}$ for some fixed element
$h\in M$ such that $\RR_{++}h \in \M$ and $S(h, h)<0$.

In Definition 1.4.1 above we gave a definition for $W$
to have arithmetic type. We can apply this definition to Weyl
group $W$ and the set of simple real roots $P(\M)$ for the Kac--Moody
algebra $\geg^{\prime\prime}(A(S,W,P(\M)))$ to have
arithmetic type. But we think that for Kac--Moody
algebras $\geg^{\prime\prime}(A(S,W,P(\M)))$ it is more
natural to give definition using imaginary roots. We will
see that these two definitions are equivalent.

The next definition is a generalization of our definition
in \cite{25} which was given for finitely generated Kac--Moody
algebras.

\definition{Definition 2.2.1} A Lorentzian Kac--Moody algebra
$\geg^{\prime\prime}(A(S,W,P(\M)))$ {\it has ari\-thmetic type}
if for any $x \in M$ with $S(x,x)<0$ there exists
$n \in \NN$ such that
$$
nx \in \Delta^{im}.
\tag2.2.1
$$
Obviously, this definition
is equivalent to the equality for open rational cones
$$
\QQ_+\Delta^{im}_+\cap V^+(S)=(M\otimes \QQ )\cap V^+(S).
\tag2.2.2
$$
\enddefinition

We have

\proclaim{Theorem 2.2.2} A Lorentzian Kac--Moody algebra
$\geg^{\prime\prime}(A(S,W,P(\M)))$ has arithmetic type if and
only if the Weyl group $W$ has arithmetic type (thus, we have
equivalent properties
\thetag{1.4.3}, \thetag{1.4.4}, \thetag{1.4.5},
\thetag{1.4.9}, \thetag{1.4.10}). Moreover, the
\newline
$\geg^{\prime\prime}(A(S,W,P(\M)))$
has arithmetic type if and only if
the Tits cone essentially coincides with $V^+(S)$ which means
$$
\text{Int}~ T=V^+(S),\ \text{equivalently,\ }\overline{T}=\overline{V^+(S)}.
\tag2.2.3
$$
\endproclaim

\demo{Proof} Assume that $\geg^{\prime\prime}(A(S,W,P(\M)))$
has arithmetic type. By \thetag{2.2.2} and
\newline
\thetag{2.1.14},
we then have
$$
(M\otimes \QQ)\cap V^+(S)\subset (\QQ_+\Delta^{im}_+)\subset T^\ast .
\tag2.2.4
$$
Since $M\otimes \QQ$ is everywhere dense in $V^+(S)$, it follows that
$V^+(S)\subset T^\ast $. Clearly, $T^\ast \subset \overline{V^+(S)}$.
It follows
the equality $\overline{T^\ast}=\overline{V^+(S)}$ and
$\overline{T}=\overline {V^+(S)}$ because $\overline{V^+(S)}^{\ \ast} =
\overline{V^+(S)}$. It follows \thetag{2.2.3}.
By \thetag{2.2.2}, \thetag{2.2.3} and \thetag{2.1.14},
we obviously get that
$(M\otimes \QQ)\cap V^+(S)\subset \QQ_+Q_+$ which is an equivalent
definition \thetag{1.4.3} for $W$ to have arithmetic type.

Now suppose that $W$ has arithmetic type. By \thetag{1.4.9} and
\thetag{1.4.10},
we then have
$$
(\RR_+Q_+)^\ast  \subset \overline{V^+(S)} \subset \RR_+Q_+ \ .
$$
It follows that $\RR_+K=(\RR_+Q_+)^\ast\cap \RR_+Q_+=
(\RR_+Q_+)^\ast \cap \overline{V^+(S)}=(\RR_+Q_+)^\ast =\RR_+\M$
where $\M$ is the fundamental polyhedron of the Weyl group $W$ in
$\La(S)$. It follows that the dual Tits cone $T^\ast=W(\RR_+K)$ contains
$V^+(S)$, and we get the equality \thetag{2.2.3}.
By \thetag{2.2.3} and \thetag{2.1.14}, we then get
$$
\QQ_+\Delta^{im}_+ \cap V^+(S) = \QQ_+ Q_+\cap V^+(S).
$$
By \thetag{1.4.3}, $V^+(S)\cap (M\otimes \QQ)\subset \QQ_+Q_+$.
Thus, by the equality above,
$V^+(S)\cap (M\otimes \QQ)\subset \QQ_+\Delta^{im}_+$. By definition,
$V^+(S) \cap (\QQ_+\Delta^{im}_+)\subset V^+(S)\cap (M\otimes \QQ)$.
Thus, we get the equality \thetag{2.2.2}, and the Lie algebra has
arithmetic type.
\enddemo

As a corollary of equivalent definitions \thetag{2.2.1},
\thetag{2.2.2} and \thetag{2.2.3},
we get

\proclaim{Corollary 2.2.3} Lorentzian Kac--Moody algebras
$\geg^{\prime\prime}(A(S,W,P(\M)))$ of arithmetic type
are exactly Lorentzian Kac--Moody algebras which have the following
two nice properties (a) and (b):

(a) Rational cone $\QQ_+\Delta^{im}_+$ generated by imaginary roots
contains $(M\otimes \QQ)\cap V^+(S)$.

(b) The denominator function $\Phi (z)$ converges in the
natural complexified cone $\Omega(S)=M\otimes \RR+iV^+(S)$ for
large  $\text{im\ }z$ (here we understand the convergence formally
as the equality of the complexified cones $\Omega (S,W)= \Omega(S)$;
for concrete cases one should support this convergence by the
appropriate estimates).

In fact, (a) and (b) are equivalent.
\endproclaim

Now let us consider Lorentzian Kac--Moody algebras of
restricted arithmetic type. One can define them like
for Definition 1.4.2.

\definition{Definition 2.2.4} A Lorentzian Kac--Moody algebra
$\geg^{\prime\prime}(A(S,W,P(\M)))$
{\it has restricted arithmetic type} if it has arithmetic type
and the symmetry group $A(P(\M))$ has finite index in $O(S)/W$
(it means that the corresponding semi-direct product
$W.A(P(\M))$ has finite index in $O(S)$).
\enddefinition

Using Theorems 1.4.3 and 2.2.2, we get

\proclaim{Theorem 2.2.5}
Let $S:M\times M\to \ZZ$ be a hyperbolic lattice, $W\subset W(S)$
an arbitrary reflection subgroup, $\M$ a fundamental polyhedron
of $W$, $P(\M)$ an acceptable set of
orthogonal vectors to $\M$ and
$$
A(P(\M))=\{ g \in O_+(S)\ | \  g(P(\M ))=P(\M ) \}
$$
the corresponding group of symmetries.
Then the Lorentzian Kac--Moody algebra
\newline
$\geg^{\prime\prime}(A(S,W,P(\M)))$ has restricted arithmetic
type if and only if the group $A(P(\M))$ of symmetries has finite index
in $O(S)/W$ (equivalently, $W.A(P(\M))$ has finite index in $O(S)$).
\endproclaim

We can rewrite the denominator formula \thetag{2.1.19} as follows:
$$
\prod_{\alpha \in \Delta_+^{im}}
{(1-e^{-2\pi i S(\alpha , z)})^{\text{mult}~(\alpha)}} =
{\sum_{w \in W} {\text{det\ }(w)e^{-2\pi i S(w(\rho )-\rho , z)}}.
\over
\prod_{\alpha \in \Delta_+^{re}}
{(1-e^{-2\pi i S(\alpha , z)})^{\text{mult}~(\alpha)}}}.
\tag2.2.5
$$
Clearly, both sides of this equality are invariant with respect to
the subgroup $W.A(P(\M))$. If the Kac--Moody
algebra $\geg^{\prime\prime}(A(S,W,P(\M)))$ has restricted
arithmetic type,
this subgroup has finite index in the arithmetic group $O(S)$.

Thus, from considerations above, we get

\proclaim{Corollary 2.2.6} Lorentzian Kac--Moody algebras
$\geg^{\prime\prime}(A(S,W,P(\M)))$ of restricted
arithmetic type are exactly Lorentzian Kac--Moody algebras which have
three nice properties (a), (b) and (c) below:

(a) Rational cone $\QQ_+\Delta^{im}_+$ generated by imaginary roots
contains $(M\otimes \QQ)\cap V^+(S)$.

(b)The denominator function $\Phi (z)$ converges in the
natural complexified cone
$\Omega(S)=M\otimes \RR+iV^+(S)$ for large $\text{im\ }z$
(here we understand the convergence formally as the
equality of the complexified cones $\Omega (S,W)= \Omega(S)$;
for concrete cases one should support this equality by the
appropriate estimates).

(c) The denominator formula \thetag{2.2.5} is invariant with respect to
the subgroup $W.A(P(\M))$ of finite index in the arithmetic group
$O(S)$.

In fact, the property (c) implies (a) and (b).
\endproclaim

At last, using results above and Theorem 1.4.7, we get that
one can always find a right sequence of extensions
\thetag{2.1.20}---\thetag{2.1.23}
of finitely generated Lorentzian Kac--Moody algebras
$\geg^{\prime\prime}(A(S, W_k, P(\M_k)))$ which improves bad properties
(-a), (-b) and (-c) of finitely generated Lorentzian
Kac--Moody algebras $\geg^{\prime\prime}(A(S, W_k, P(\M_k)))$.

\proclaim{Theorem 2.2.7}
Let $S:M\times M\to \ZZ$ be a hyperbolic lattice, $W_1\subset W(S)$
a finitely generated reflection subgroup, $\M_1$ a fundamental polyhedron
of $W_1$ with a finite acceptable set $P(\M_1)$ of vectors orthogonal
to $\M_1$. Let $\geg^{\prime\prime}(A(S, W_1, P(\M_1)))$ be the
corresponding Lorentzian Kac--Moody algebra.
Then $W_1$ is a reflection subgroup $W_1\subset W \subset W(S)$ of
a reflection group $W$ of restricted arithmetic type with a
fundamental polyhedron $\M\subset \M_1$ and an acceptable set
$P(\M)$ such that $P(\M_1)\subset P(\M)$. This gives
the embedding
$$
\geg^{\prime\prime}(A(S, W_1, P(\M_1)))\subset
\geg^{\prime\prime}(A(S, W, P(\M)))
$$
of the finitely generated Lorentzian Kac--Moody algebra
$\geg^{\prime\prime}(A(S, W_1, P(\M_1)))$
to the Lorentzian Kac--Moody algebra
$\geg^{\prime\prime}(A(S, W, P(\M)))$
of restricted arithmetic type (which has
nice properties (a), (b) and (c) above).
\endproclaim

\subhead
2.3. Lorentzian Kac--Moody algebras of
restricted arithmetic type with a lattice Weyl vector
\endsubhead

Results above show that Lorentzian Kac--Moody algebras of
restricted arithmetic type give the most natural class of
Lorentzian Kac--Moody algebras to study. But this class is
huge. Every hyperbolic lattice $S$ of rank $\ge 3$
with not-trivial reflection group $W(S)$  has
infinitely many reflection subgroups $W\subset W(S)$ of
restricted arithmetic type, which have infinite
index between one another. Each this group defines its own Lorentzian
Kac--Moody algebras of restricted arithmetic type
$\geg^{\prime\prime}(A(S, W, P(\M)))$.

Fortunately, there is another very natural restriction on
Lorentzian Kac--Moody algebras.

The denominator
formula has the best and the most interesting
automorphic properties only if it can be written in
the form \thetag{2.1.15} when one can consider real and imaginary roots
together. This may be done only if there exists a lattice
Weyl vector $\rho$.

\definition{Definition 2.3.1} Lorentzian Kac--Moody algebra
$\geg^{\prime\prime}(A(S, W, P(\M)))$ has  {\it a lattice Weyl vector} if
there exists $\rho \in M\otimes \QQ$ such that
$$
S(\rho, \alpha)=-(1/2)S(\alpha, \alpha)\ \text{for any\ } \alpha \in P(\M).
\tag2.3.1
$$
\enddefinition

If the lattice Weyl vector $\rho$ does exist, we can write the denominator
formula in the form \thetag{2.1.15} which we rewrite here
$$
\Phi (z)=e^{-2\pi iS(\rho,z)} \prod_{\alpha \in \Delta_+}
{(1-e^{-2\pi i S(\alpha , z)})^{\text{mult}~(\alpha)}} =
\sum_{w \in W} {\text{det\ }(w)e^{-2\pi i S(w(\rho ), z)}}.
\tag2.3.2
$$
It is anti-invariant with respect to $W.A(P(\M))$ which means that
$$
\Phi (w.a(z))=det(w)\Phi (z)\ \text{for any\ } w\in W,a\in A(P(\M)).
\tag2.3.3
$$

Using results above, we get

\proclaim{Theorem 2.3.2}
Lorentzian Kac--Moody algebras
$\geg^{\prime\prime}(A(S,W,P(\M)))$ of restricted
arithmetic type and with lattice Weyl vector
are exactly Lorentzian Kac--Moody algebras which have
three nice properties (a), (b) and (c)+(d) below:

(a) Rational cone $\QQ_+\Delta^{im}_+$ generated by imaginary roots
contains $(M\otimes \QQ)\cap V^+(S)$.

(b)The denominator function $\Phi (z)$ converges in the complexified cone
$\Omega(S)=M\otimes \RR+iV^+(S)$ for big $\text{im\ }z$
(here we understand the convergence formally as the
equality of the complexified cones $\Omega (S,W)= \Omega(S)$;
for concrete cases one should support this equality by the
appropriate estimates).

(c)+(d) The denominator formula \thetag{2.3.2}
is anti-invariant with respect to
the subgroup $W.A(P(\M))$ of finite index in the arithmetic group
$O(S)$.
\endproclaim

We apply results of Sect. 1.1 and Sect. 1.4 to describe these the most
interesting Lorentzian Kac--Moody algebras.

Like in Sect. 1.1, we have

\definition{Definition 2.3.2}
Let $S:M\times M\to \ZZ$ be a hyperbolic lattice, $W\subset W(S)$
a reflection subgroup, $\M$ a fundamental polyhedron
of $W$, $P(\M)$ an acceptable set of
orthogonal vectors to $\M$ and
$$
A(P(\M))=\{ g \in O_+(S)\ | \  g(P(\M ))=P(\M ) \}
$$
the corresponding group of symmetries.
Lorentzian Kac--Moody algebra \linebreak
$\geg^{\prime\prime}(A(S,W,P(\M)))$ is {\it elliptic} if
$W\subset O(S)$ has finite index. In particular, the lattice
$S$ is {\it elliptic reflective}.
Lorentzian Kac--Moody algebra
$\geg^{\prime\prime}(A(S,W,P(\M)))$
is {\it parabolic} if
$A(P(\M))$ is infinite but $W.A(P(\M))$ has finite index in
$O(S)$ and there exists $c \in M$ such that $S(c,c)=0$ and
$g(c)=c$ for any $g \in A(P(\M))$ (i.e.  $A(P(\M))$ has
a cusp $c$). In particular, the lattice $S$
is {\it parabolic reflective}.
\enddefinition

By Theorems 1.3.2 and 1.3.3, Propositions 1.4.10 and 1.4.11 and
results of Sect. 2.2, we get

\proclaim{Theorem 2.3.3}
Lorentzian Kac--Moody algebras
$\geg^{\prime\prime}(A(S,W,P(\M)))$
of restricted arithmetic type and with a lattice Weyl vector
$\rho$ have $S(\rho, \rho)\le 0$ if $\rk S\ge 3$.

For $S(\rho, \rho)<0$, these algebras are
exactly elliptic Lorentzian Kac--Moody algebras
with a lattice Weyl vector. Their set is finite.

For $S(\rho, \rho)=0$, these algebras are
exactly parabolic Lorentzian Kac--Moody algebras
with a lattice Weyl vector. For a constant $D>0$, their set is finite if
$$
[O(S)_\rho:A(P(\M))]<D.
$$
\endproclaim

\vskip10pt

It seems that existence of the generalized lattice Weyl vector $\rho$
(see Definition 1.4.9) is
also very important for automorphic properties of the denominator
formula \thetag{2.3.2} of Lorentzian
Kac--Moody algebras. (The denominator formula is well-defined after
a finite-dimensional extension of $M$ by the subspace in the kernel of
the Gram matrix of $P(\M)$.)

By Theorems 1.1.1 and 1.1.3, Propositions 1.4.10 and 1.4.11,
and results of Sect. 1.3.2, we get

\proclaim{Theorem 2.3.4}
Lorentzian Kac--Moody algebras
$\geg^{\prime\prime}(A(S,W,P(\M)))$
of restricted arithmetic type and with a generalized
lattice Weyl vector $\rho$  have $S(\rho, \rho)\le 0$ if $\rk S\ge 3$.

For $S(\rho, \rho)<0$, these algebras are
exactly elliptic Lorentzian Kac--Moody algebras.
The set of their lattices $S:M\times M \to \ZZ$ is finite.

For $S(\rho, \rho)=0$, these algebras are
exactly elliptic or parabolic Lorentzian Kac--Moody algebras.
The set of their lattices $S:M\times M\to \ZZ$ is also finite.
\endproclaim

\remark{Remark 2.3.5}
These results show that classification of elliptic and parabolic
reflective lattices is of the extremal importance for the theory
of Lorentzian Kac--Moody algebras. Many elliptic reflective
lattices were constructed by \'E.B. Vinberg, see \cite{29} ---
\cite{34}. In \cite{16}, \cite{17},
\cite{18}, \cite{22} all elliptic 2-reflective even hyperbolic lattices
$S$ were classified. Here $S$ is elliptic $2$-reflective if the subgroup
$W^{(2)}(S)$ generated by reflections in all $\delta \in M$ with
$S(\delta, \delta)=2$ has finite index in $O(S)$.
In \cite{17} there is a series of twelve even 2-elementary
(i.e. with $2$-elementary discriminant group $M^\ast /M$) parabolic
2-reflective even hyperbolic lattices.
The even unimodular hyperbolic lattice of the
rank $26$ which is parabolic by J. Conway \cite{7}
continues this series.
In R. Borcherds \cite{1} there are some examples of
elliptic and parabolic reflective lattices. See reviews of
\'E.B. Vinberg, O.V. Shvartsman \cite{34} and
R. Scharlau, C. Walhorn \cite{28}
for further results on reflective lattices.

We should say that in spite of finiteness results of \S 1 above,
complete description of reflective hyperbolic lattices is the very
difficult problem. On the other hand, now there are plenty of reflective
lattices known, and because of these finiteness results,
all these examples are very interesting.
\endremark

\subhead
2.4. A comment on the recent results due to R. Borcherds
\endsubhead

Here we want to write down how we understand the recent results of
R. Borcherds \cite{2} --- \cite{6}
in connection with the results above.

Let us consider a Lorentzian Kac--Moody algebra
$\geg^{\prime\prime}(A(S,W,P(\M)))$ of restricted arithmetic type
and with a lattice Weyl vector $\rho$ (thus, by Sect. 2.3, this
Lorentzian Lie algebra is elliptic or parabolic).
Then its denominator function $\Phi (z)$ is defined in the
complexified cone $V^+(S)$ which is
$$
\Omega (S)=M\otimes \RR+iV^+(S),
\tag2.4.1
$$
and converges for $\im z>>0$.

We denote by $U:N\times N\to \ZZ$ an even unimodular
lattice of the signature $(1,1)$ where
$N=\ZZ e_1\oplus \ZZ e_2$ and
$U(e_1,e_1)=U(e_2,e_2)=0,\ U(e_1,e_2)=-1$.

We consider a lattice
$$
S^\prime =S\oplus U(k):M^\prime\times M^\prime \to \ZZ,\ \ k \in \NN,
\tag2.4.2
$$
where $M^\prime =M\oplus N$.
The lattice $S^\prime$ defines the type IV domain  which is one
of two connected components of the domain
$$
\Omega (S^\prime)=\{\CC\omega \subset M^\prime \otimes \CC\
|\ S^\prime (\omega, \omega )=0,\
S^\prime(\omega , \overline\omega )<0 \}.
\tag2.4.3
$$
We normalize $\omega\in \CC\omega \in \Omega(S^\prime)$ by the
condition $S^\prime (\omega , e_1)=-1$. Then
$$
\omega =z\oplus (S(z,z)/2)e_1 \oplus (1/k)e_2,\ z \in \Omega (S).
\tag2.4.4
$$
This defines an embedding
$$
\Omega (S) \subset \Omega (S^\prime)
$$
which is called {\it of the cusp $e_1$ embedding}.

We comment the recent results by R. Borcherds as follows:

\vskip5pt

{\it The most interesting Lorentzian Kac--Moody
algebras $\geg^{\prime\prime}(A(S,W,P(\M)))$
of restricted arithmetic type and with a lattice
Weyl vector $\rho$ (thus, they are elliptic or parabolic) have
a generalized Kac--Moody algebra (or may be generalized
Kac--Moody superalgebra) correction
$\geg^{\prime\prime}(A(S,W,P(\M)))\subset
\geg^{\prime\prime}(A(S,W,P(\M)\cup H))$ such
that the denominator function $\Phi^\prime (z)$ of the
correcting generalized Kac--Moody algebra
$\geg^{\prime\prime}(A(S,W,P(\M)\cup H))$
is an automorphic form
with respect to an appropriate subgroup of finite index in
$O(S^\prime)$ for some $k$.} It is natural to name this correction as
{\it automorphic form correction} (or {\it automorphic correction}).

\vskip5pt

We remark that a priory the denominator function $\Phi (z)$ is
anti-invariant with respect to the subgroup of finite index
$M^\ast_\rho.(W.A(P(\M))$ of the stabilizer subgroup
$O(S^\prime)_{e_1}$ of the cusp. Here $M^\ast_\rho=
\{x \in M^\ast|S(\rho, x)\in \ZZ\}$.

More precisely this means that there exists an automorphic form
$\widetilde{\Phi}(z)$ which is automorphic with
respect to a subgroup of finite index $G\subset O(S^\prime)^+$,
containing $W.A(P(\M))$, and is anti-invariant with respect to
$W.A(P(\M))$. Moreover, $\widetilde{\Phi}(z)$
has Fourier expansion with the properties which we describe below.
This Fourier expansion is
$$
\widetilde{\Phi}(z)=
\sum_{w\in W}
{\det (w)\left( e^{-2\pi iS(w(\rho ),z)}
-\sum_{a\in M^\ast_{P(\M)}\cap \RR_{++}\M}
{m(a)e^{-2\pi iS(w(\rho+a),z)}}\right)} ,
\tag2.4.5
$$
where all $m(a)\in \ZZ$ and $m(a) \in \ZZ_+$
(for superalgebras case $m(a)\in \ZZ$) if $S(a,a)<0$.
See \thetag{1.2.4} for definition of $M^\ast_{P(\M)}$.
Let us consider $a \in M^\ast_{P(\M)}\cap \RR_{++}\M$
such that $S(a,a)=0$. These $a$ define several isotropic rays
which correspond to infinite vertices of $\M$,
and the cusp of $\M$ if $\M$ is parabolic. We consider any of
these rays $ta_0$, where $a_0$ is a primitive element
of $M^\ast_{P(M)}$ with $S(a_0,a_0)=0$ and $t \in \NN$.
For any of these isotropic rays one has the equality of formal
power series of one variable $q$
$$
1-\sum_{t \in \NN}{m(ta_0)q^t}=
\prod_{k\in \NN}{(1-q^k)^{\tau(ka_0)}},
\tag2.4.6
$$
where $\tau (ka_0) \in \ZZ_+$
(for superalgebras case $\tau(ka_0)\in \ZZ$).

Let us consider
$$
\split
H=&\{m(a)a\ |\ a \in M^\ast_{P(\M)} \cap \RR_{++}\ \M \ \text{and}\
S(a,a)<0\}\cup\\
\cup &\{\tau(a)a\ |\ a \in M^\ast_{P(\M)} \cap \RR_{++}\ \M \
\text{and}\ S(a,a)=0 \},
\endsplit
\tag2.4.7
$$
where $ka$ means that we repeat an element $a$ exactly $k$ times
(consider with the multiplicity $k$).
We consider the Gram matrix
$$
G(P(\M)\cup H)
\tag2.4.8
$$
of elements $G(P(\M)\cup H)$. One can easily see that
$G(P(\M)\cup H)$ has all necessary properties to define
a generalized Kac--Moody algebra
$$
\geg^{\prime\prime}(G(S,W,P(\M)\cup H))
\tag2.4.9
$$
corresponding to this matrix (see R. Borcherds \cite{2} for
definition).
It has properties similar to Kac--Moody algebras we
considered in Sects 2.1 and 2.2.
The only difference is that its
simple roots are $P(\M)\cup H$.
It is graded by the lattice $M^\ast_{P(\M)}\supset M$, has roots,
multiplicities of roots and the same Weyl group $W$. The Kac--Moody
algebra $\geg^{\prime\prime}(A(S,W,P(\M)))$ is defined
by Gram matrix $G(P(\M))$ of elements $P(\M)$ and is
a subalgebra of $\geg^{\prime\prime}(G(S,W,P(\M)\cup H))$.
The denominator formula for $\geg^{\prime\prime}(G(S,W,P(\M)\cup H))$
equals to
$$
e^{-2\pi iS(\rho,z)} \prod_{\alpha \in \Delta_+}
{(1-e^{-2\pi i S(\alpha , z)})^{\text{mult}~(\alpha)}} =
\widetilde{\Phi}(z)
\tag2.4.10
$$
where $\widetilde{\Phi} (z)$ is the function described in
\thetag{2.4.5}---\thetag{2.4.7} (see
\cite{3} where in fact necessary calculations were done).
In particular, if one can find the automorphic
function $\widetilde{\Phi}(z)$ which we described above,
it automatically has the product formula
\thetag{2.4.10}.

The next the most amazing example of the "correction" above
was found by R. Borcherds \cite{3}, \cite{5}.
Consider an even unimodular hyperbolic
lattice $S$ of signature $(25,1)$
and $W=W(S)$. This reflection group is parabolic and has the
lattice Weyl vector $\rho$ with $S(\rho, \rho )=0$. This
was proved by J. Conway \cite{7}.
Take $P(\M)=P(\M)_{pr}$ for a fundamental polyhedron $\M$ of $W$.
The corrected denominator formula for this algebra equals
$$
\split
\widetilde{\Phi}(z)=&e^{-2\pi iS(\rho,z)} \prod_{\alpha \in \Delta_+}
{(1-e^{-2\pi i S(\alpha , z)})^{p_{24}(1-S(\alpha,\alpha)/2)}} = \\
&\sum_{w\in W}{\det(w)\sum_{n>0}{\tau (n)e^{-2\pi in(w(\rho), z)}}} .
\endsplit
\tag2.4.11
$$
Here
$$
\sum_{n\ge 0}{p_{24}(n)q^n}=\prod_{n>0}{(1-q^n)^{-24}}
\tag2.4.12
$$
and
$$
\sum_{n\ge 0}{\tau (n)q^n}=q\prod_{n>0}{(1-q^n)^{24}}.
\tag2.4.13
$$
Thus, for this case $H=\rho, 2\rho,....,t\rho,...$,
where each element $t\rho$ is taken with multiplicity $24$.
The function $\widetilde{\Phi}(z)$ is an automorphic form of weight $12$
with respect to the subgroup  $O(S\oplus U)^+\subset O(S\oplus U)$
which fixes connected components of the domain $\Omega(S\oplus U)$.
See some other examples in \cite{5}, \cite{6}
(and also \cite{3}, \cite{4}) and \cite{12}.

It seems that the bases for a construction of similar
examples and many others is the arithmetic lifting of
Jacobi forms on IV type domains which was
constructed by V.A. Gritsenko \cite{10}, \cite{11}. Some its
multiplicative analog was constructed by R. Borcherds \cite{6}.

\vskip10pt

Theory of automorphic forms is a very delicate domain.
Automorphic forms like \thetag{2.4.11} are probably very rare.
On the other hand, as we have shown above,
essentially one has only a finite number
of elliptic and parabolic Kac--Moody algebras with
a lattice Weyl vector. Thus, essentially one
has to "correct" only a finite set of Kac--Moody
algebras. May be it is possible like for the example above.

\enddocument

\end